\documentclass{amsart}
\usepackage{amssymb} % enables mathbb etc.
\usepackage{mathtools} % enables xleftarrow etc.
\usepackage{float} %for proper working of hyperref with floats like algorithm; must appear before hyperref.
\usepackage{hyperref}
\usepackage{cleveref}
\usepackage{algorithm} %float wrapper for algorithms, not compatible with hyperref without float, hyperref, algorithm order.
\usepackage{enumitem}
\usepackage{array}
\usepackage[noend]{algorithmic} %algorithm type setting environment
\numberwithin{algorithm}{subsection}
\numberwithin{equation}{subsection}
\numberwithin{table}{subsection}

\theoremstyle{plain}
\newtheorem{conjecture}{Conjecture}[subsection]
\newtheorem{lemma}[conjecture]{Lemma}
\newtheorem{theorem}[conjecture]{Theorem}
\newtheorem{corollary}[conjecture]{Corollary}
\newtheorem{proposition}[conjecture]{Proposition}

\theoremstyle{remark}
\newtheorem{remark}[conjecture]{Remark}

\theoremstyle{definition}
\newtheorem{definition}[conjecture]{Definition}

\DeclareMathOperator{\trd}{trd}
\DeclareMathOperator{\End}{End}
\DeclareMathOperator{\Gal}{Gal}

\DeclareMathOperator{\nrd}{nrd}

\DeclareMathOperator{\Idl}{Idl}
\DeclareMathOperator{\PIdl}{PIdl}
\DeclarePairedDelimiterX{\norm}[1]{\lVert}{\rVert}{#1}
\DeclarePairedDelimiterX{\round}[1]{\lfloor}{\rceil}{#1}
\newcommand{\curlyO}{\mathcal{O}}
\newcommand{\Fp}{\mathbb{F}_p}
\newcommand{\Fpp}{\mathbb{F}_{p^2}}
\newcommand{\mi}{\mathbf{i}}
\newcommand{\mj}{\mathbf{j}}
\newcommand{\mk}{\mathbf{k}}

\makeatletter
\let\@wraptoccontribs\wraptoccontribs
\makeatother

\title{Gross lattices of supersingular elliptic curves}

\author{Chenfeng He}
\address{E\"{o}tv\"{o}s Lor\'{a}nd University, Budapest, Hungary}
\email{chenfenghe163@gmail.com}

\author{Gaurish Korpal}
\address{University of Arizona, Tucson, United States}
\curraddr{University of Auckland, Auckland, New Zealand}
\email{gkorpal@arizona.edu}

\author{Ha T. N. Tran}
\address{University of Alberta -- Augustana Campus, Camrose, Canada}
\email{htran2@ualberta.ca}

\author{Christelle Vincent}
\address{University of Vermont, Burlington, United States}
\email{christelle.vincent@uvm.edu }

\contrib[with an appendix by]{Jonathan Love}

\thanks{ The authors thank Jonathan Love for stimulating questions that led to the statement and proof of \Cref{prop:D_3_1728} and \Cref{thm:OqvsOqprime}, as well as the referees of an earlier version of the article for their insightful questions and comments. This work started as a collaboration of the authors with Victoria de Quehen and Sarah Days-Merrill during the Isogeny Graphs in Cryptography Workshop 2023 at the Banff International Research Station for Mathematical Innovation and Discovery (BIRS). In addition, this article is based on work supported by the National Science Foundation under Grant No. DMS-1928930, while the authors were in residence at the Simons Laufer Mathematical Sciences Institute in Berkeley, California, during the Summer Research in Mathematics program 2024. Ha T. N. Tran acknowledges the support of the Natural Sciences and Engineering Research Council of Canada (NSERC) (funding RGPIN-2019-04209 and DGECR-2019-00428). Finally, Christelle Vincent acknowledges the hospitality of the GAATI laboratory at the University of French Polynesia during part of her work on this article.}

\subjclass[2020]{Primary 11G20, 11R52, 14G15, 14G50;
Secondary 11H06, 11Y40, 11Y16.}
\begin{document}

\begin{abstract}
Let $p$ be a prime, $E$ be a supersingular elliptic curve defined over $\overline{\mathbb{F}}_p$, and $\curlyO$ be its (geometric) endomorphism ring. Earlier results of Chevyrev-Galbraith and Goren-Love have shown that the successive minima of the Gross lattice of $\curlyO$ characterize the isomorphism class of $\curlyO$. In this paper, we extend this work and show that the value of the third successive minimum $D_3$ of the Gross lattice gives necessary and sufficient conditions for the curve to have its $j$-invariant in the field $\Fp$ or in the set $\Fpp \setminus \mathbb{F}_p$, as well as finer information about the endomorphism ring of $E$ when its $j$-invariant belongs to $\mathbb{F}_p$ and $p \equiv 3 \pmod{4}$. We end our article with an investigation of the geometry of Gross lattices of supersingular elliptic curves.
\end{abstract}

\maketitle
%%%%%%%%%%%%%%%%%%%%%%%%%

\section{Introduction}\label{sec:intro}

Let $p$ be a prime, $B_p$ be the quaternion algebra over $\mathbb{Q}$ that ramifies exactly at $p$ and $\infty$, and $E$ be a supersingular elliptic curve defined over $\overline{\mathbb{F}}_p$. Then the geometric endomorphism ring $\End_{\overline{\mathbb{F}}_p}(E)$ of $E$ is known to be isomorphic to a maximal order $\curlyO$ in $B_p$, and furthermore, every such isomorphism class of maximal orders corresponds to the endomorphism ring of a unique $\Gal(\overline{\mathbb{F}}_p/\mathbb{F}_p)$-orbit of isomorphism classes of supersingular elliptic curves over $\overline{\mathbb{F}}_p$. The original proof of this correspondence was given by Deuring in \cite{DeuringCorrespondence}, and a particularly compact and complete account of what we now accordingly call the \emph{Deuring correspondence} appears in Table 1 of \cite{SQISign}.

Interestingly, this correspondence can only be computed efficiently in one direction: Given an isomorphism class of maximal orders in $B_p$ we know how to compute a model for a supersingular elliptic curve defined over $\overline{\mathbb{F}}_p$ with endomorphism ring isomorphic to this order in time polynomial in $\log(p)$ \cite{wesolowski1,deuring}. Conversely, however, the best algorithms to determine the isomorphism class of the endomorphism ring of a supersingular elliptic curve as an order in $B_p$ solely from an equation for $E$ have exponential running time in $\log(p)$ \cite{EHLMP2,FIKMN}. Since the introduction of so-called \emph{supersingular elliptic curve cryptography} in \cite{CharlesLauterGoren} and then \cite{JaoDeFeo,JaoDeFeoPlut}, the implications of the existence of this correspondence and the fact of its ``one-wayness" on the security of this family of cryptographic systems has been an important, as well as mathematically quite technical, object of study  \cite{EHLMP,Bank2019,BonehLove,wesolowski2,PW24,merdy2023supersingular}. 

\subsection{Main results of this article}

A classical tool to obtain information about the isomorphism class of a maximal order $\curlyO$ in $B_p$ is its so-called \emph{Gross lattice} $\curlyO^T$, a rank-$3$ sublattice containing the trace zero elements of $\curlyO$ of the form $2x- \trd(x)$ for $x$ in $\curlyO$, and which was first defined by Gross \cite[Proposition 12.9]{Gross}. Throughout, if $E$ is a supersingular curve defined over $\overline{\mathbb{F}}_p$ with endomorphism ring isomorphic to $\curlyO$, we will say that $\curlyO^T$ is the Gross lattice of $E$; the relationship of a supersingular elliptic curve to its Gross lattice has been the object of several articles since, and we briefly highlight the results we will need from these works in Section \ref{sec:previously}. Much of the work we draw on has developed in particular results on the values $D_1 \leq D_2 \leq D_3$ commonly called the \emph{successive minima} of the Gross lattice $\curlyO^T$, defined here in \Cref{def:succmin}. In this article we continue the investigation of the Gross lattices of maximal orders in the quaternion algebra $B_p$, and especially of their geometry, by which we mean not only their successive minima but also finer invariants, such as their \emph{shape}, in the sense of \cite{Terr} and \cite{BhargavaHarris}.

We begin in Section \ref{sec:mainsection} by extracting information about the field of moduli of a supersingular elliptic curve from the successive minima of its Gross lattice, and more specifically, the value of the third successive minimum $D_3$. Our first main result can be obtained by combining Theorems \ref{thm:D3notoverFp} and \ref{thm:D3overFppart1}, as well as Proposition \ref{prop:0}: 

\begin{theorem}\label{thm:main1}
Let $E$ be a supersingular elliptic curve defined over $\overline{\mathbb{F}}_p$, $\curlyO^T$ be its Gross lattice, and $D_3$ be the third successive minimum of $\curlyO^T$. If $p \geq 7$, then $j(E)$, the $j$-invariant of $E$, belongs to $\mathbb{F}_p\setminus\{0\}$ if and only if 
\begin{equation*}
p \le D_3 \le \frac{8p}{7} + \frac{7}{4};
\end{equation*}
and if $p \neq 3$, $j(E) = 0$ if and only if $D_3= \frac{4p+1}{3}$.
Otherwise, $j(E)$ belongs to $\mathbb{F}_{p^2} \setminus \mathbb{F}_p$ and 
\begin{equation*}D_3 \le \frac{3p}{5} + 5.\end{equation*}
\end{theorem}

The statement of this theorem invites two quick remarks before we move on: First, in the case of $p = 3$, there is a unique supersingular elliptic curve defined over $\overline{\mathbb{F}}_3$, its $j$-invariant is equal to $0 \equiv 1728 \pmod{3}$ and belongs to $\mathbb{F}_3$, and in this case $D_3 = 4$ (see Remark \ref{rem:1728}). Secondly, one might wonder about the cases where $p < \frac{3p}{5} + 5$, which creates an overlap between the two bounds stated above. In this case, we have $p \leq 11$, and every supersingular elliptic curve defined over $\overline{\mathbb{F}}_p$ has $j$-invariant in $\mathbb{F}_p$ satisfying $D_3 \geq p$.
To demonstrate the strength of this result, we additionally show that it is very close to best possible in Theorem \ref{thm:OqvsOqprime} and Propositions~\ref{prop:Fp2_D1is20} and \ref{prop:jminus15cube} .

The second main result of Section \ref{sec:mainsection} shows that the third successive minimum $D_3$ yields finer information about $\curlyO$ when it is the endomorphism ring of a supersingular elliptic curve with $j$-invariant in $\mathbb{F}_p$ and $p \equiv 3 \pmod{4}$.
Indeed, we characterize the third succesive minimum of the curves whose endomorphism ring contains an element of norm $\frac{p+1}{4}$ in Theorem~\ref{thm:OqvsOqprime}; this result was also obtained independently in \cite[Proposition 3]{Clements} using different methods:
\begin{theorem}\label{thm:main2}
    Let $p \equiv 3 \pmod{4}$, $E$ be a supersingular elliptic curve defined over $\overline{\mathbb{F}}_p$ with $j$-invariant in $\mathbb{F}_p$, $\curlyO^T$ be its Gross lattice, and $D_3$ be the third successive minimum of $\curlyO^T$. Then $\curlyO$, the endomorphism ring of $E$, is maximally embedded by $\mathbb{Z}[\frac{1+\sqrt{-p}}{2}]$ if and only if $D_3 \in \{p,p+1\}$. This latter case, $D_3= p+1$, occurs if and only if $j(E) = 1728$.
\end{theorem}

We note here that our results in Theorems \ref{thm:main1} and \ref{thm:main2} thus suggest that supersingular elliptic curves with $j$-invariant in $\mathbb{F}_p$ (whose corresponding vertices are often called the \emph{spine} of the supersingular isogeny graph, after \cite{supersingularland}) have more small-degree endomorphisms than those with $j$-invariants in $\mathbb{F}_{p^2} \setminus \mathbb{F}_p$, which is a result independently obtained by Orvis \cite[Corollary 4.2]{Orvis}. As remarked by Love in a recent presentation on the topic \cite{lovetalk}, it may be interesting to see if our results may be obtained from those of Orvis or vice-versa.
 
We now turn to presenting the second set of main results of this article, which pertain to the geometry, or shape, of the Gross lattice of a supersingular elliptic curve, and are contained in Section \ref{sec:geometry}. We begin in Section \ref{subsec:mainorthorounded} by investigating if the Gross lattice of a supersingular elliptic curve can be orthogonal or well-rounded. By establishing necessary and sufficient conditions for the successive minima of the Gross lattice of a supersingular elliptic curve to be equal and for vectors of a \emph{successive minimal basis} (see Definition \ref{def:succminbasis}) to be orthogonal, we find that this is generally impossible, except when $p=2$ in which case the Gross lattice of the unique isomorphism class of supersingular elliptic curves defined over $\overline{\mathbb{F}}_2$ is well-rounded. 

Following this, we turn our attention in Section \ref{sec:uniqueness} to the so-called \emph{normalized Gram matrix} of a supersingular elliptic curve, which is the Gram matrix of a normalized successive minimal basis for its Gross lattice (both introduced in Definition \ref{def:normalized}) and is of the form given in \Cref{lem:Gram_matrix}. Our interest in this quantity comes from the fact that we can show in \Cref{thm:Gram_unique}, relying on a result of \cite{Goren-Love}, that in the case where $p \neq 3$, the Gram matrices of all normalized successive minimal bases of $\curlyO^T$ are equal and thus, this matrix is a well-defined invariant for isomorphism classes of supersingular elliptic curves when $p \neq 3$.

Finally, at this end of this article we turn our attention to the case of supersingular elliptic curves with $j$-invariant in $\mathbb{F}_p$. In this setting, we show in Section~\ref{sec:algorithm} that the normalized Gram matrix can only of one of four types (given in Theorem \ref{thm:4types}); we present in the appendix a proof of this result that was suggested to us by Jonathan Love as it is much shorter and simpler than the one we had obtained. This theorem allows us to give a simple algorithm to compute the possible normalized Gram matrices of a supersingular elliptic curve with $j$-invariant in $\mathbb{F}_p$ given $p$, the characteristic of its base field, and $D_1$, the value of the first successive minimum of its Gross lattice. 

In Section~\ref{sec:special_cuves}, the last section of this article, we apply this algorithm in  to compute the normalized Gram matrix of $\curlyO^T$ for certain special families of supersingular elliptic curves. Indeed, there are precisely $13$ so-called \emph{CM $j$-invariants} -- values $j \in \overline{\mathbb{Q}}$ corresponding to elliptic curves with complex multiplication -- that in fact belong to $\mathbb{Q}$ \cite{heegner,stark}. Accordingly, each value corresponds to an elliptic curve whose $\overline{\mathbb{Q}}$-isomorphism class contains a curve defined over $\mathbb{Q}$, and with complex multiplication by one of the $13$ imaginary quadratic orders of class number $1$. Upon reducing modulo $p$ for suitable primes  $p$, these curves yield supersingular elliptic curves with $j$-invariants in $\mathbb{F}_p$, and we establish a lower bound on $p$ that ensures that the generator of the endomorphism ring in characteristic zero corresponds to the shortest vector in the Gross lattice of its reduction modulo $p$. In each case, the algorithm yields a unique possible Gram matrix, which we conclude is the Gram matrix of a normalized successive minimal basis for the Gross lattice of the elliptic curve.

\subsection{Computational tools and structure of the article}
Some of the results below rely on computational software to handle the cases of small primes. To calculate $\mathbb{Z}$-bases for maximal orders in $B_p$, we use the algorithm given by Kirschmer and Voight~\cite{KV10} and available in \texttt{Magma}~\cite{magma}. To compute the successive minima of Gross lattices, we apply an implementation of Eisenstein reduction~\cite[Theorem 103]{Dickson}  written by Gustavo Rama and available within \texttt{SageMath}\footnote{\url{https://github.com/sagemath/sage/blob/develop/src/sage/quadratic_forms/ternary_qf.py}}. All further computations of Gross lattices and their Gram matrices were done using our own functions written using \texttt{SageMath}~\cite{sagemath}, and our code is available for review on GitHub at~\cite{repo}.

The structure of this article is as follows: We begin in \Cref{sec:backgound} by introducing the notation and definitions we will need and establishing and recalling some results we will need on successive minimal bases and Gross lattices. Following this, our main results on the third successive minimum of Gross lattices are presented in \Cref{sec:mainsection}. Finally, in \Cref{sec:geometry}, we delve into the study of the geometry of the Gross lattice.

%%%%%%%%%%%%%%%%%%%%%%%%%%%

\section{Background}\label{sec:backgound}

We begin by introducing the facts, notations, and definitions we will need in this article.

\subsection{Quaternion algebras}
Throughout this article, $p$ will always be a prime, and recall that $B_p$ is the quaternion algebra over $\mathbb{Q}$ that ramifies exactly at $p$ and $\infty$. Then $B_p$ is equipped with a standard involution $\overline{\cdot}$ which we call \textbf{conjugation}, and we denote the \textbf{reduced trace} of $x \in B_p$ by $\trd(x) = x+\overline{x}$ and the \textbf{reduced norm} of $x \in B_p$ by $\nrd(x) = x\overline{x}$; the reduced norm is a positive definite quadratic form on $B_p$. In addition, $B_p$ is equipped with an \textbf{inner product} $( x, y ) = \frac{1}{2} \trd(x \overline{y})$. The norm $\norm{\cdot}$ associated to this inner product is the square root of the usual reduced norm $\nrd$ on $B_p$; indeed since $\nrd(x) = x \overline{x} \in \mathbb{Q}$, $\trd(x \overline{x}) = 2 \nrd(x)$ so $\norm{x}^2=(x, x ) = \nrd(x)$. In this article, we will privilege the notation $\norm{x}^2$ over $\nrd(x)$, and call this quantity the \textbf{norm} of $x$.

An element $x$ of $B_p$ is \textbf{integral} if and only if both $\trd(x)$ and $\norm{x}^2$ are in $\mathbb{Z}$. The set of integral elements of $B_p$ in general does not form a ring, and the algebra $B_p$ can contain more than one  maximal orders of integral elements.

Since every element of $B_p \setminus \mathbb{Q}$ generates an imaginary quadratic field over $\mathbb{Q}$, one way to study the maximal orders of $B_p$ is to identify the imaginary quadratic orders $O$ that can be embedded in a given maximal order $\mathcal{O}$ of $B_p$, and more specifically the quadratic orders that can be \emph{optimally embedded} in $\mathcal{O}$:

\begin{definition}\label{def:maximally_embedded}
    Let $O$ be an order in an imaginary quadratic field $K$ over $\mathbb{Q}$. We say that a maximal order $\curlyO$ of $B_p$ is \textbf{maximally embedded} by $O$ if there exists an embedding $\iota \colon K \to B_p$ such that  $\iota(O) =\curlyO \cap \iota(K)$. Such an embedding is called an \textbf{optimal embedding} of $O$ in $\curlyO$.
\end{definition}

We also take this opportunity to give a proof of the following result, which as far as we could find was first proved in \cite[Proposition 69]{kohelthesis} using different techniques, and which is a slight generalization of \cite[Lemma 4.1]{Chevyrev-Galbraith}. It is also a variation on, and generalization of, \cite[Lemma 3.9]{Goren-Love}, which is proved using the same ideas we use here:

\begin{lemma}\label{lem:unique_order}
Let $O$ be an imaginary quadratic order of class number $1$. If two maximal orders $\curlyO$ and $\curlyO'$ of $B_{p}$ are maximally embedded by $O$, then $\curlyO$ and $\curlyO'$ are isomorphic.
\end{lemma}
\begin{proof}
Throughout this proof we use the notation and results developed in \cite{JVbook}. Let $K = O \otimes_{\mathbb{Z}} \mathbb{Q}$ be the quadratic field containing $O$, and for $p$ a prime we write
\begin{equation*}
\left( \frac{K}{p} \right) = 
\begin{cases}
1 & \text{if $p$ splits in $K$;}\\
0 & \text{if $p$ ramifies in $K$;}\\
-1 & \text{if $p$ is inert in $K$.}
\end{cases}
\end{equation*}
In addition, for $\curlyO$ a maximal order in $B_p$, we denote by $m(O,\curlyO,\curlyO^\times)$ the number of optimal embeddings of $O$ into $\curlyO$, up to conjugation by elements of $\curlyO^\times$. Finally for $S$ either $O$ or $\curlyO$, we denote by $h(S)$ the cardinality of the quotient $\Idl(S)/\PIdl(S)$, consisting of the invertible fractional ideals of $S$ modulo its principal two-sided ideals ($\Idl(S)$ is abelian in both cases considered here).

Then putting together equations (30.4.21) and (30.7.2) and Proposition 30.5.3(b) of \cite{JVbook}, and specializing them to the case at hand, we obtain the formula
\begin{equation}\label{eq:countembeds}
\sum_{\curlyO \text{ maximal}} h(\curlyO) m(O,\curlyO,\curlyO^\times) = h(O) \left( 1- \left( \frac{K}{p} \right) \right),
\end{equation}
and every term of the equation is a nonnegative integer. Of course, this formula recovers the fact that if $\left( \frac{K}{p} \right)=1$, then there is no embedding of $K$ into $B_p$, and so we do not consider this case here.

If $p$ ramifies in $K$ and $h(O) = 1$, the right hand side of equation \eqref{eq:countembeds} is equal to $1$, and therefore only one term in the sum on the left can be positive. Hence only one isomorphism class of maximal order in $B_p$ maximally embeds $O$.

If $p$ is inert in $K$ and $h(O) = 1$, the right hand side of equation \eqref{eq:countembeds} is equal to $2$, and we prove that in this case every term in the sum on the left is even, which again implies that only one of them can be positive, and only one isomorphism class of maximal order in $B_p$ maximally embeds $O$.

If $\curlyO^\times = \{\pm 1\}$, then $m(O,\curlyO,\curlyO^\times)$ simply counts the number of optimal embeddings of $O$ into $\curlyO$, as conjugation by $\curlyO^\times$ is trivial. Since embeddings of $O$ into $\curlyO$ comes in conjugate pairs, and the complex conjugate of an optimal embedding is also optimal, $m(O,\curlyO,\curlyO^\times)$ is even. This proves the result for every $p \not \equiv 11 \pmod{12}$, since when we avoid this case there is at most one maximal order in $B_p$ with nontrivial units, and the term corresponding to this order must also be even since every other term in the equation is even. To handle the case of $p \equiv 11 \pmod{12}$, when two maximal orders of $B_p$ have nontrivial units, we note that if $\# \curlyO^\times = 6$, then computing the eigenvalues of the action by conjugation by elements of $\curlyO^\times$ on $B_p$, we see that none of the units act by multiplication by $-1$ on any subspace of $B_p$, and hence conjugation by a unit cannot identify an optimal embedding with its complex conjugate, and again $m(O,\curlyO,\curlyO^\times)$ is even. This completes the proof.
\end{proof}
\begin{remark}
The proof of Lemma \ref{lem:unique_order} shows that if $p \equiv 3 \pmod{4}$ is inert in an imaginary quadratic field $K$ containing $O$ an order of class number $1$ and $\# \curlyO^\times = 4$, then $h(\curlyO)m(O,\curlyO,\curlyO^\times)$ is again even. Note that this does not contradict the fact that $h(\curlyO)m(O,\curlyO,\curlyO^\times)$ must be odd for some maximal order $\curlyO$ when $p$ is ramified; indeed if $p \geq 5$ ramifies in $K$ and $K$ contains an order $O$ of class number $1$, then $p \equiv 3 \pmod{4}$, $B_p$ contains a unique maximal order with  $\# \curlyO^\times = 4$, this order is maximally embedded by $O$, and the action by conjugation of the elements of $\curlyO^\times$ do induce complex conjugation on the image of this embedding, yielding $m(O,\curlyO,\curlyO^\times) = 1$ in this case.
\end{remark}

\subsection{The Gross lattice}
Given $\curlyO$ a maximal order in $B_p$, the main object of study in this article is its Gross lattice:
\begin{definition}\label{def:grosslattice}
Let $\curlyO$ be a maximal order in $B_p$. Its  \textbf{Gross lattice} $\curlyO^T$ is the sublattice given by    
\begin{equation}\label{eq:grosslattice}
\curlyO^T = \{ 2x - \trd(x) : x \in \curlyO\}.
\end{equation}
When $E$ is a supersingular elliptic curve and $\End(E) \cong \curlyO$, we also call $\curlyO^T$ the \textbf{Gross lattice of $E$}.
\end{definition}

The Gross lattice $\mathcal{O}^T$ of a maximal order in $B_p$ is a free $\mathbb{Z}$-module of rank $3$, and throughout we write $\{\beta_1, \beta_2, \beta_3\}$ for an ordered $\mathbb{Z}$-basis of $\curlyO^T$. 
We note in addition that the bilinear pairing $(\gamma_1,\gamma_2) = \frac{1}{2} \trd(\gamma_1 \overline{\gamma}_2)$ takes integer values for any two elements $\gamma_1,\gamma_2$ belonging to a given Gross lattice. Indeed, if $\gamma_i = 2 \alpha_i - \trd(\alpha_i)$ for $\alpha_i$ in a maximal order $\curlyO$, then 
\begin{equation*}
\frac{1}{2} \trd(\gamma_1 \overline{\gamma}_2) = 2 \trd(\alpha_1 \overline{\alpha}_2)- \trd(\alpha_1) \trd(\alpha_2),
\end{equation*}
which belongs to $\mathbb{Z}$ since $\mathcal{O}$ is closed under conjugation and thus $\alpha_1 \overline{\alpha}_2$ is integral if $\alpha_1$ and $\alpha_2$ belong to the same maximal order.

The significance of the Gross lattice of a maximal order $\mathcal{O}$ in $B_p$ is the following:

\begin{proposition}[{\cite[Proposition 3.6]{Goren-Love}}] \label{prop:endo_degd}
Let $p$ be a prime and $\mathcal{O}$ be a maximal order in $B_p$. Then there is an embedding of $O_{-d}$, the imaginary quadratic order of discriminant $-d$, in $\curlyO$ if and only if its Gross lattice $\curlyO^T$ contains an element of norm $d$. Furthermore, $\curlyO$ is maximally embedded by $O_{-d}$ if and only if the corresponding element of norm $d$ is a primitive element of $\curlyO^T$, that is, it is not a nontrivial integer multiple of any other element of $\curlyO^T$.
\end{proposition}

As a consequence of Proposition \ref{prop:endo_degd}, we note that if $\beta \in \mathcal{O}^T$, then $\norm{\beta}^2$ is congruent to $0$ or $3$ modulo $4$.
This is a fact which we will use again and again throughout this article. 

\subsection{Lattice definitions and notations}

Throughout this article, a lattice $\Lambda$ in $B_p$ with ordered $\mathbb{Z}$-basis $\{b_1, \hdots, b_n\}$ is denoted by $\langle b_1, \hdots, b_n\rangle$. Furthermore, we define the \textbf{Gram matrix} of this basis to be the symmetric matrix 
\begin{equation*}
G_{\{b_1, \hdots, b_n\}}= \left((b_i, b_j)\right)_{i, j} = \left(\frac{1}{2}\trd(b_i \overline{b_j})\right)_{i, j},
\end{equation*}
and the \textbf{determinant} of $\Lambda$, denoted by $\det(\Lambda)$, can be given by the quantity $\det(G_{\{b_1, \hdots, b_n\}})$ (which is independent of the choice of basis). 

\begin{definition}\label{def:succmin}
Let $\Lambda$ be a lattice of rank $n$  in $B_p$. For $1 \leq i \leq n$, we define the \textbf{$i$th successive minimum} of $\Lambda$ to be the smallest value $D_i$ such that the rank of the $\mathbb{Z}$-submodule of $\Lambda$ generated by $\{x \in \Lambda : \norm{x}^2 \leq D_i\}$ is greater than or equal to $i$. 
\end{definition}
We note that the above definition does not agree with the standard definition of successive minima from lattice theory, as we use the quantity $\norm{x}^2$ rather than $\norm{x}$ to define them, following \cite{Chevyrev-Galbraith} and \cite{Goren-Love}. 

\begin{definition}\label{def:succminbasis}
Let $\Lambda$ be a lattice of rank $n$ and $D_1 \leq \cdots \leq D_n$ its successive minima. An ordered list of elements $\{x_1, \ldots, x_n\} \in \Lambda$ \textbf{attains the successive minima of $\Lambda$} if $\norm{x_i}^2 = D_i$ for each $i$. A lattice $\Lambda$ of rank at most $3$ always has a basis that attains its successive minima \cite[Corollary 2.6.10]{martinet13}, and we call such an ordered basis a \textbf{successive minimal basis} of $\Lambda$. 
\end{definition}

We note that by \cite[Corollary 6.2.3]{martinet13} in fact any list of $3$ elements attaining the successive minima of a lattice of rank $3$ is a basis for this lattice, and this result is also true for any order in $B_p$ if $p$ is odd by \cite[Lemma 3.5]{Goren-Love}.

Now, for a lattice $\Lambda$ of rank $3$ in $B_p$, which will be the main case of interest in this article, we may compute a successive minimal basis by computing an Eisenstein-reduced basis for $\Lambda$ (which gives us a reduced  fundamental parallelopiped for $\Lambda$ as in \cite[pp.\ 162-163]{Dickson}). By the construction of the Eisenstein-reduced basis given in \emph{ibid.}, such a basis is automatically Minkowski-reduced \cite{Minkowski}, and since we are in dimension smaller than $4$, a Minkowski-reduced basis attains the successive minima of $\Lambda$ \cite{vanderWaerden}.

Using the notation above, standard lattice bounds show that if $\Lambda$ is a lattice of rank $n$, there is a minimal constant $\gamma_n$ (called the $n$-th Hermite constant) such that
\begin{equation}\label{eq:hermite}
    \det(\Lambda) \leq \prod_{i=1}^n D_i \leq \gamma_n^n \det(\Lambda)
\end{equation}
(\cite[Theorem 2.6.8]{martinet13} gives the upper bound) and we have $\gamma_2^2 = \frac{4}{3}$ and $\gamma_3^3 = 2$.

To end this section we recall the size-reducedness condition for vectors in a basis of a lattice, and present an important lemma on the size-reducedness of pairs of elements in successive minimal bases for lattices of rank $3$:

\begin{definition}\label{def:size_reduced}
    Given an ordered basis $\{b_1, b_2, \cdots, b_n\}$ of a lattice $\Lambda$, we can apply the Gram–Schmidt process to obtain an orthogonal basis $\{b_1, b_2^*, \cdots, b_n^*\}$ for $\Lambda \otimes_\mathbb{Z} \mathbb{R}$, which we call the \textbf{Gram-Schmidt orthogonalization} of $\{b_1, b_2, \cdots, b_n\}$. Furthermore given this ordered basis $\{b_1, b_2, \cdots, b_n\}$ we define its \textbf{Gram-Schmidt coefficients} to be
    \begin{equation*}
    \mu_{j,i} = \frac{(b_j, b_i^*)}{(b_i^*, b_i^*)} \quad \text{for} \quad i < j.
    \end{equation*}
    Finally, the ordered pair $\{b_i, b_j\}$ for $i<j$ is called \textbf{size-reduced} if $|\mu_{j,i}| \le \frac{1}{2}$.
\end{definition}
Note that if a pair $\{b_i, b_j\}$ is not size-reduced, then we can obtain a new size-reduced pair $\{b_i, b'_j\}$ by replacing $b_j$ with $b'_j = b_j -\lfloor \mu_{j,i} \rceil b_i$, where $\lfloor \mu_{j,i} \rceil$ denotes the integer closest to $\mu_{j,i}$ and the value $0.5$ is rounded down to $0$. Moreover, $\langle b_i, b_j\rangle = \langle b_i, b'_j\rangle$, i.e., both pairs generate the same lattice. For more details, see \cite{LLL82}.
 
\begin{lemma}\label{lem:size_reduced}
If $\Lambda'$ is a lattice of rank $2$ with $\{v_1, v_2\}$ a successive minimal basis, then the pair $\{v_1, v_2\}$ is size-reduced. Similarly, if $\Lambda$ is a lattice of rank $3$ with $\{v_1, v_2, v_3\}$ a successive minimal basis, then the pairs $\{v_1, v_2\}$ and  $\{v_1, v_3\}$  are size-reduced. Moreover, in this case, $$\left|\frac{(v_3, v_2)}{(v_2, v_2)}\right| \le \frac{1}{2}.$$  
\end{lemma}

\begin{proof}
Throughout we use the notation introduced in \Cref{def:size_reduced}, and note that it suffices to prove the results for $\Lambda$ a lattice of rank $3$ since the size-reducedness of the pair $\{v_1, v_2\}$ attaining the first two successive minima of a lattice does not depend on the dimension of the ambient lattice.

First, suppose on the contrary that $\{v_1, v_2\}$ is not size-reduced, so $|\mu_{2,1}| > \frac{1}{2}$. Let $v'_2= v_2- \lfloor \mu_{2,1} \rceil v_1 \in \Lambda$, then we have
\begin{align*}
\norm{v'_2}^2 & = \|v_2\|^2 + (\lfloor \mu_{2,1} \rceil)^2\|v_1\|^2 - 2 \lfloor \mu_{2,1} \rceil (v_2, v_1) \\
&= \|v_2\|^2 + \lfloor \mu_{2,1} \rceil (\lfloor \mu_{2,1} \rceil-2 \mu_{2,1}) \|v_1\|^2.
\end{align*}
Furthermore
\begin{equation*}
| \lfloor \mu_{2,1} \rceil | - |\mu_{2,1}| \leq \frac{1}{2} < |\mu_{2,1}|,
\end{equation*}
and considering the cases of $\mu_{2,1}$ positive and negative separately, this implies that $\lfloor \mu_{2,1} \rceil (\lfloor \mu_{2,1} \rceil-2 \mu_{2,1})<0$. Therefore, $\|v'_2\|^2 < \|v_2\|^2$ and the subset $\{v_1, v'_2\}$ of $\Lambda$ is linearly independent. This contradicts the fact that $\|v_2\|^2$ is the second successive minimum of $\Lambda$.

Turning our attention to the last two statements of the lemma, assuming that $|\mu_{3,1}| > \frac{1}{2}$ and setting $v'_3= v_3- \lfloor \mu_{3,1} \rceil v_1 \in \Lambda$, we can show that $\|v'_3\|^2 < \|v_3\|^2$, and the subset $\{v_1, v_2, v'_3\}$ of $\Lambda$ is linearly independent, similarly to the case worked out above. These contradict the fact that $\|v_3\|^2$ is the third successive minimum of $\Lambda$. Finally, if $\delta= \frac{(v_3, v_2)}{(v_2, v_2)}$ and $ |\delta| > \frac{1}{2}$, we can show that $\|v'_3\|^2 < \|v_3\|^2$ where $v'_3= v_3- \lfloor \delta \rceil v_2 \in \Lambda$ and $\{v_1, v_2, v'_3\}$ is linearly independent, which again is a contradiction.
\end{proof}

\subsection{Supersingular elliptic curves}
Throughout, if a variety or morphism can be described by an equation with coefficients in a field $k$, we say that this variety or morphism is \textbf{defined over} that field $k$. 
Furthermore, $E$ will most often be a supersingular elliptic curve defined over $\overline{\mathbb{F}}_p$ for $p$ a prime. Such an elliptic curve necessarily has $j$-invariant $j(E)$ in $\Fpp$, and an element of its $\overline{\mathbb{F}}_p$-isomorphism class is defined over $\mathbb{F}_{p^2}$; and when $j(E)$ belongs to $\Fp$, then there is a curve in the $\overline{\mathbb{F}}_p$-isomorphism class of $E$ that is defined over $\Fp$. The set of elements of $\mathbb{F}_{p^2}$ that are $j$-invariants of supersingular elliptic curves are called \textbf{supersingular $j$-invariants}. 

If $E$ is a supersingular elliptic curve, then the \textbf{geometric endomorphism ring} of $E$, which contains the endomorphisms of $E$ defined over $\overline{\mathbb{F}}_p$ and which we denote $\End(E)$ throughout, is isomorphic to a maximal order $\curlyO$ in the quaternion algebra $B_p$ over $\mathbb{Q}$ ramified at $p$ and $\infty$. Conveniently, if $\phi \in \End(E)$ is mapped to $\alpha \in \curlyO$ under an isomorphism $\End(E) \cong \curlyO$, then the degree of $\phi$ is equal to $\nrd(\alpha)$, and its trace is equal to $\trd(\alpha)$.  

For any supersingular $j$-invariant in $\mathbb{F}_p$ and $E$ a supersingular elliptic curve defined over $\mathbb{F}_p$ with this $j$-invariant, the \textbf{arithmetic endomorphism ring} of $E$ is the subring of its geometric endomorphisms that are defined over $\mathbb{F}_p$. Interestingly, except if $j(E) = 1728$, the isomorphism class of the arithmetic endomorphism ring of a supersingular elliptic curve defined over $\mathbb{F}_p$ is determined by the $j$-invariant of $E$, and not its $\mathbb{F}_p$-isomorphism class. Since a variety defined over $\mathbb{F}_p$ must admit a $p$-power Frobenius endomorphism, which is an endomorphism of degree $p$ and trace $0$, the arithmetic endomorphism ring of an elliptic curve defined over $\mathbb{F}_p$ must admit an embedding by $\mathbb{Z}[\sqrt{-p}]$, and therefore is an order $O$ in the quadratic field $\mathbb{Q}(\sqrt{-p})$. We note that this order is not necessarily maximal if $p \equiv 3 \pmod{4}$, as in this case $\frac{1+\sqrt{-p}}{2}$ is also integral. In any case, if $O$ is the arithmetic endomorphism ring of $E$ and $\curlyO$ its geometric endomorphism ring, then $\curlyO$ is maximally embedded by $O$. 

Finally, if $E$ is an elliptic curve defined over $\overline{\mathbb{Q}}$ with endomorphism ring isomorphic to an imaginary quadratic order $O$ -- in which case we say that $E$ has \textbf{complex multiplication} by $O$ -- and $\overline{\mathfrak{P}}$ is a prime of $\overline{\mathbb{Q}}$ of good reduction for $E$, then the endomorphism ring of the reduction of $E$ modulo $\overline{\mathfrak{P}}$ is maximally embedded by $O$. In particular, if $K = O \otimes_\mathbb{Z} \mathbb{Q}$ is the imaginary quadratic field containing $O$ and $p = \overline{\mathfrak{P}} \cap \mathbb{Z}$ is inert or ramified in $K$, then the reduction of $E$ modulo $\overline{\mathfrak{P}}$ is supersingular, say with endomorphism ring isomorphic to $\curlyO$ in $B_p$, and $\curlyO$ is maximally embedded by $O$. (If the reduction of $E$ modulo $\overline{\mathfrak{P}}$ is ordinary, then its endomorphism ring is isomorphic to $O$.)

\subsection{Previous work on the topic}\label{sec:previously}

Finally, for the convenience of the reader, in this section we recall and sometimes restate various results that we will use. Throughout, for $E$ a supersingular elliptic curve, $D_1 \leq D_2 \leq D_3$ are the successive minima of its Gross lattice (we recall that if $\End(E) \cong \mathcal{O}$ for $\mathcal{O}$ a maximal order in $B_p$, for simplicity we call $\mathcal{O}^T$ the Gross lattice of $E$).

We begin with a fundamental result on the structure of the endomorphism ring of supersingular elliptic curves with $j(E) \in \mathbb{F}_p$ due to Ibukiyama:

\begin{theorem}[{\cite[Theorems 1 and 2, and Lemmata 1.2 and 1.8]{Ibukiyama}}]\label{thm:ibukiyama}
Let $p$ be an odd prime and $B_p$ be the quaternion algebra over $\mathbb{Q}$ ramified at $p$ and $\infty$. Let $q \equiv 3 \pmod{8}$ be a prime such that $\left(\frac{-q}{p}\right) = -1$ and $r$ be an integer such that $r^2 + p \equiv 0 \pmod{q}.$ Then $B_p$ is isomorphic to $\mathbb{Q}+\mi \mathbb{Q} + \mj\mathbb{Q} + \mi\mj\mathbb{Q}$ where $\mi^2 = -p, \mj^2 = -q$ and $\mi\mj = -\mj\mi$, and the order
\begin{equation}
\curlyO(q,r) = \mathbb{Z} + \frac{1+\mj}{2}\mathbb{Z} + \frac{\mi(1+\mj)}{2}\mathbb{Z} + \frac{(r+\mi)\mj}{q}\mathbb{Z}
\end{equation}
is a maximal order in $B_p$, whose isomorphism class does not depend on the choice of $r$, and which optimally embeds the imaginary quadratic order $\mathbb{Z}[\sqrt{-p}]$. If $p \equiv 1 \pmod{4}$, every maximal order in $B_p$ is isomorphic to an order of the form $\curlyO(q,r)$ as $q$ varies. 

If $p \equiv 3 \pmod{4}$, for $q' \equiv 3 \pmod{8}$ a prime such that $\left(\frac{-q'}{p}\right) = -1$ and $r'$ an integer such that $(r')^2 +p \equiv 0 \pmod{4q}$, taking again $\mi$ such that $\mi^2=-p$ and $\mj$ such that $\mj^2 = -q'$ with $\mi \mj = -\mj\mi$, then the order
\begin{equation}
\curlyO'(q',r') = \mathbb{Z} + \frac{1+\mi}{2}\mathbb{Z} + \mj\mathbb{Z} + \frac{(r'+\mi)\mj}{2q'}\mathbb{Z}
\end{equation}
is also maximal, its isomorphism class does not depend on the choice of $r'$, and $\curlyO'(q',r')$ optimally embeds the quadratic order $\mathbb{Z}[\frac{1+\sqrt{-p}}{2}]$. If $p \equiv 3 \pmod{4}$, every maximal order in $B_p$ is isomorphic to an order of the form $\curlyO(q,r)$ or $\curlyO'(q',r')$ as $q$ and $q'$ vary. Furthermore, an order of the form $\curlyO(q,r)$ for some $q$ and $r$ is isomorphic to one of the form $\curlyO'(q',r')$ for some $q'$ and $r'$ if and only if both $\curlyO(q,r)$ and $\curlyO'(q',r')$ contain four units. 
\end{theorem}

We also note this result of Elkies, which implies that $D_1 < p$ if $p \geq 11$:

\begin{proposition}[{\cite[Section 4]{Elkies}}]\label{prop:elkies}
Let $p$ be a prime and $E$ any supersingular elliptic curve defined over $\overline{\mathbb{F}}_p$. Then $D_1$, the first successive minimum of its Gross lattice, satisfies
\begin{equation}
D_1 \leq 2p^{2/3}.
\end{equation}
\end{proposition}

When $j(E) \in \mathbb{F}_p$, Kaneko obtains the following refinement of this bound:

\begin{theorem}[{\cite[Theorem 1]{kaneko1989}}]\label{thm:kaneko}
Let $p$ be a prime and $E$ any supersingular elliptic curve defined over $\overline{\mathbb{F}}_p$ with $j(E) \in \mathbb{F}_p$. Then $D_1$, the first successive minimum of its Gross lattice, satisfies
\begin{equation}
D_1 \leq \frac{4}{\sqrt{3}}\sqrt{p}.
\end{equation} 
\end{theorem}

\begin{remark}
    \label{rem:uniquej}
    Since every supersingular curve defined over $\overline{\mathbb{F}}_p$ has $j(E) \in \mathbb{F}_p$ when $p \leq 31$, this refined bound shows that $D_1 < p$ for $p \geq 7$. We note in addition that both results above are best possible by \cite[Theorem 1.2 and Proposition 1.4]{YangMinCM}. 
\end{remark}

Returning to Kaneko's work, in the same article \cite{kaneko1989} the following two results are also implicitly shown, which are crucial for this work:

\begin{proposition}[discussed in {\cite[Section 3]{Chevyrev-Galbraith}}, stated here as in {\cite[Proposition 3.12]{Goren-Love}}] 
    Let $p$ be a prime and $\mathcal{O}$ be an order in $B_p$ (in particular $\mathcal{O}$ need not be maximal), and let $\gamma_1, \gamma_2 \in \curlyO^T$ be linearly independent. Then 
\begin{equation*}
\norm{\gamma_1}^2\norm{\gamma_2}^2 - \frac{1}{4}\trd(\gamma_1 \overline{\gamma_2})^2
\end{equation*}
is a positive integer multiple of $4p$. As a consequence, every rank-$2$ sublattice of $\curlyO^T$ has determinant $4np$ for some positive integer $n$.
\label{prop:sublattice4np}
\end{proposition} 

As well as: 

\begin{proposition}[implicit in {\cite[pages 851-852]{kaneko1989}}, stated here as in the proof of {\cite[Lemma 2.3]{Chevyrev-Galbraith}}] 
Let $p$ be a prime and $E$ any supersingular elliptic curve defined over $\overline{\mathbb{F}}_p$ with $j(E) \in \mathbb{F}_p$. Then its Gross lattice contains a sublattice of rank 2 and of determinant $4p$.  
\label{prop:sublat4pkaneko}
\end{proposition}

\begin{proof}
Since we did not find the proof explicitly anywhere, we give it here: if $E$ has endomorphism ring isomorphic to an order of the form $\curlyO(q,r)$, then the lattice with basis $\gamma_1 = \mj, \gamma_2 = \frac{2(r+\mi)\mj}{q}$, where one can verify that $\gamma_1,\gamma_2 \in \mathcal{O}^T$, has determinant $4p$, and if $E$ has endomorphism ring isomorphic to an order of the form $\curlyO'(q',r')$, then the lattice with basis $\gamma_1 = 2\mj, \gamma_2 = \frac{(r'+\mi)\mj}{q'}$ has determinant $4p$, and again one can verify that $\gamma_1,\gamma_2 \in \mathcal{O}^T$ using the explicit basis for $\curlyO'(q',r')$.
\end{proof}

Another important result for us is the following:

\begin{lemma}[{\cite[Lemma 3.1]{Chevyrev-Galbraith}}]
    Let $\mathcal{O}$ be a maximal order in $B_p$, then $\det(\mathcal{O}^T) = 4p^2$.
    \label{lem:detOT}
\end{lemma}

Using the Hermite bound of \eqref{eq:hermite}, this yields the following inequality, where as usual $D_1 \leq D_2 \leq D_3$ are the successive minima of the Gross lattice of a maximal order in $B_p$:
\begin{equation}\label{eq:prodbound}
4p^2 \leq D_1 D_2 D_3 \leq 8p^2.
\end{equation}

We also note this result, whose converse we will show in Remark \ref{rem:CGconverse}:

\begin{lemma}[{\cite[Lemma 2.3]{Chevyrev-Galbraith}}]\label{lem:CGd1d2}
Let $\mathcal{O}$ be a maximal order in $B_p$, $D_1,$ $D_2$ the first two successive minima of $\mathcal{O}^T$, and $E$ a supersingular elliptic curve with endomorphism ring isomorphic to $\mathcal{O}$. If $j(E) \in \mathbb{F}_p$, then 
\begin{equation}
D_1D_2 \leq \frac{16p}{3}.
\end{equation}
\end{lemma}

By Remark \ref{rem:CGconverse}, we have that \cite{Chevyrev-Galbraith} proves the following exactly when $j(E) \in \mathbb{F}_p$. The full result, due to \cite{Goren-Love}, shows the significance of the successive minima of the Gross lattice of a supersingular elliptic curve, and motivated the work of Section \ref{sec:mainsection}:

\begin{theorem}[{\cite[Theorem 2.1]{Chevyrev-Galbraith}} for the case of $j(E) \in \mathbb{F}_p$, the full result is a corollary of {\cite[Theorem 1.4]{Goren-Love}}]
Let $\mathcal{O}$ and $\mathcal{O}'$ be two maximal orders in $B_p$. If their Gross lattices $\mathcal{O}^T$ and $\mathcal{O}'^T$ have the same successive minima, then $\mathcal{O}$ and $\mathcal{O}'$ are isomorphic.
\end{theorem}

To prove \cite[Theorem 1.4]{Goren-Love}, the authors obtain the following two results; they suffice for their purposes. As a consequence of our work, we can remove certain hypotheses to state them more generally. First, in Proposition \ref{prop:inequalityDi} we slightly improve the following result to remove the hypotheses that $p \geq 11$ and $D_1 \geq 15$, and replace them with $p \geq 3$ and $D_1 \geq 4$: 

\begin{theorem}[adapted from {\cite[Lemmata 4.4 and 4.5]{Goren-Love}}]\label{thm:GLinequality}
Let $p \geq 11$ be a prime, $E$ be a supersingular elliptic curve defined over $\overline{\mathbb{F}}_p$ whose Gross lattice has successive minima $D_1 \leq D_2 \leq D_3$. Then $D_1 \neq D_2$. Furthermore, if $D_1 \geq 15$, then $D_2 \neq D_3$ as well.   
\end{theorem}

Finally, we expand on the following, removing the condition that $D_1 \geq 8$, suggesting a normalization for the values $\frac{1}{2}\trd(\beta_i \overline{\beta}_j)$ that yields a unique normalized Gram matrix except when $p=3$, and improving the bound on $\left| \frac{1}{2}\trd(\beta_i \overline{\beta}_j)\right|$. These results can be found in Lemma \ref{lem:Gram_matrix}, Proposition \ref{prop:tijzero} and Theorem \ref{thm:Gram_unique}.

\begin{theorem}[Corollary 3.15 as well as Section 4.1 of \cite{Goren-Love}]
 Let $p$ be an odd prime, $\mathcal{O}$ be a maximal order in $B_p,$ and let $\{\beta_1,\beta_2,\beta_3\}$ be a successive minimal basis of $\mathcal{O}^T$. If $\min \{ \norm{\beta_i}^2,\norm{\beta_j}^2 \} \leq p,$ then the absolute value of the bilinear product $T_{ij} =\left| \frac{1}{2}\trd(\beta_i \overline{\beta}_j)\right|$ is the unique integer square root modulo $p$ of $\norm{\beta_i}^2 \norm{\beta_j}^2$ in the interval $[0,\frac{p}{2}]$. 
 
 As a consequence, if now $\mathcal{O}^T$ has first successive minimum $D_1 \geq 8$, then the Gram matrix of any successive minimal basis is equal to the Gram matrix of any other, up to the sign changes induced on the values $\frac{1}{2}\trd(\beta_i \overline{\beta}_j)$ resulting from sending $\beta_i$ to $\pm\beta_i$ for $1\leq i \leq 3$.
 \label{thm:GLgrammatrix}
\end{theorem}

\subsection{The Gram matrix of a successive minimal basis}

As announced, we end this section with a result which strengthens Theorem \ref{thm:GLgrammatrix} slightly, by improving the bound on the absolute value of $(\beta_i,\beta_j) = \frac{1}{2}\trd(\beta_i \overline{\beta}_j)$ when $\beta_i$ and $\beta_j$ are elements of a successive minimal basis, and by fixing the values of $(\beta_1,\beta_2)$ and $(\beta_1,\beta_3)$ to be nonnegative (rather than $(\beta_1,\beta_2)$ and $(\beta_2,\beta_3)$, as in \cite{Goren-Love}). The reason for this choice will be made clear as a result of Proposition \ref{prop:tijzero}, as explained in Remark \ref{rem:choice}. We note that the improved bound on the absolute value of $(\beta_i,\beta_j)$ can also be obtained from the proof of \cite[Corollary 3.15]{Goren-Love}.

\begin{lemma}\label{lem:Gram_matrix}
Let $E$ be a supersingular curve defined over $\overline{\mathbb{F}}_p$ with Gross lattice $\curlyO^T$ and let $\{\beta_1, \beta_2, \beta_3\}$ be a successive minimal basis of $\curlyO^T$. By possibly replacing $\beta_2$ with $-\beta_2$ and $\beta_3$ with $-\beta_3$ to ensure that $(\beta_1,\beta_2)$ and $(\beta_1,\beta_3)$ are nonnegative, we obtain that the Gram matrix $G_{\{\beta_1, \beta_2, \beta_3\}}$ of this (new) basis can be written in the form 
\begin{equation}\label{eq:Gram_mattrix}    G_{\{\beta_1, \beta_2, \beta_3\}}=\begin{pmatrix}
D_1 & t_{12} & t_{13} \\
t_{12} & D_2 & t_{23} \\
t_{13} & t_{23} & D_3
\end{pmatrix},
\end{equation}
where $D_i= \|\beta_i\|^2$, $t_{ij}=(\beta_i, \beta_j)\in \mathbb{Z}$,  $ 0\le t_{12}, t_{13}  \le \frac{D_1}{2}$ and $|t_{23}| \le\frac{D_2}{2}$. 
Furthermore, if $p$ is odd and $D_1 \geq 8$, imposing the condition that $0\le t_{12}, t_{13}$ defines a unique Gram matrix $G_{\curlyO^T}$ which is independent of the choice of basis.
\end{lemma}
 \begin{proof}   
This follows from \Cref{lem:size_reduced} since $t_{12}= \mu_{2,1}D_1, t_{13}= \mu_{3,1}D_1$ and $t_{23}= \frac{(\beta_3, \beta_2)}{(\beta_2, \beta_2)} D_2$, and because the inner product of two elements in a given Gross lattice is an integer.
Finally, the uniqueness follows from Theorem \ref{thm:GLgrammatrix}.
\end{proof}

This result as well as Proposition \ref{prop:tijzero}, as explained in Remark \ref{rem:choice}, prompts the following definition:

\begin{definition}\label{def:normalized}
We say that a successive minimal basis $\{\beta_1, \beta_2, \beta_3\}$ for $\curlyO^T$ is \textbf{normalized} if the inner products $(\beta_1,\beta_2)$ and $(\beta_1,\beta_3)$ are nonnegative. The Gram matrix of a normalized successive minimal basis is a \textbf{normalized Gram matrix} for $\mathcal{O}^T$.
\end{definition}

%%%%%%%%%%%%%%%%%%%%%%%%%%%%%%%%%%%%%%%%%%%%%%%%
\section{Successive minima of Gross lattices}\label{sec:mainsection}

This section contains the results characterizing the field of definition of the $j$-invariant of a supersingular elliptic curve in terms of the third successive minimum $D_3$ of its Gross lattice. After presenting some results on supersingular elliptic curves whose endomorphism ring is optimally embedded by an imaginary quadratic order of class number $1$ in Section \ref{sec:class1}, we show in \Cref{subsec:rank2sub} that a necessary and sufficient condition for a supersingular elliptic curve to have $j$-invariant in $\Fp$ is that its Gross lattice $\curlyO^T$ has a rank-$2$ sublattice of determinant $4p$, and in addition, we show that such a sublattice can be generated by two elements that achieve the first two successive minima of $\curlyO^T$, a result crucial to the work that follows. 
Following this, \Cref{subsec:0and1728} contains the computation of successive minimal bases as well as their Gram matrices for curves with $j$-invariants $0$ and $1728$ when they are supersingular. With these preliminaries in hand, we can finally show our main results on the third successive minimum of $\curlyO^T$ for each case: The case of $j(E) \notin \Fp$ is handled in \Cref{subsec:D3notinFp} and the case of $j(E) \in \Fp$ in \Cref{subsec:D3inFp}. 
\subsection{Maximal orders optimally embedded by imaginary quadratic orders of class number $1$}\label{sec:class1}

We will need the following consequence of Lemma \ref{lem:unique_order}:

\begin{corollary}\label{cor:ssreduction}
Let $O$ be an imaginary quadratic order of class number $1$, and $E$ be a supersingular elliptic curve defined over $\overline{\mathbb{F}}_p$ such that its endomorphism ring is maximally embedded by $O$. Then $j(E)$, the $j$-invariant of $E$, belongs to $\mathbb{F}_p$.
\end{corollary}

\begin{proof}
We first note that by equation \eqref{eq:countembeds}, $p$ must be inert or ramified in $K = O \otimes_{\mathbb{Z}} \mathbb{Q}$, the imaginary quadratic field containing $O$.

Now, let $\alpha$ be such that $O = \mathbb{Z}[\alpha]$ and let $\phi$ be an endomorphism of $E$ that is the image of $\alpha$ under an optimal embedding $O \hookrightarrow \End(E)$. By Deuring's lifting theorem \cite{DeuringCorrespondence}, or \cite[1.7.4.5]{ChaiConradOort} for a more precise statement, the pair $(E, \phi)$ lifts to characteristic zero, by which we mean that there is an elliptic curve $\widetilde{E}$ defined over $\overline{\mathbb{Q}}$ which has good reduction at a prime $\overline{\mathfrak{P}}$ of $\overline{\mathbb{Q}}$ with $p = \overline{\mathfrak{P}} \cap \mathbb{Z}$ and such that the reduction modulo $\overline{\mathfrak{P}}$ of $\widetilde{E}$ is $\overline{\mathbb{F}}_p$-isomorphic to $E$, and the elliptic curve $\widetilde{E}$ has an endomorphism $\tilde{\phi}$ such that the reduction modulo $\overline{\mathfrak{P}}$ of $\tilde{\phi}$ is $\phi$. Therefore, in particular, the imaginary quadratic order $O' \cong \End_{\overline{\mathbb{Q}}}(\widetilde{E})$ contains $O$.

As a consequence, $O'$ has class number $1$, and $j(\widetilde{E}) \in \mathbb{Z}$. Since $j(E) \equiv j(\widetilde{E}) \pmod{p}$, we conclude that $j(E) \in \mathbb{F}_p$.
\end{proof}

A simple corollary of this result along with Proposition \ref{prop:endo_degd} is the following, which we will need repeatedly throughout this article:

\begin{corollary}\label{cor:0and1728}
    Let $E$ be a supersingular elliptic curve defined over $\overline{\mathbb{F}}_p$. If its Gross lattice $\curlyO^T$ contains an element $\beta$ with $\norm{\beta}^2 = 3$, then $j(E) = 0$, and if its Gross lattice $\curlyO^T$ contains an element $\beta$ with $\norm{\beta}^2 = 4$, then $j(E) = 1728$.
\end{corollary}

\begin{proof}
In either case, if there is an element $\beta$ in $\curlyO^T$ of norm equal to $3$, or $4$ respectively, $\beta$ must be primitive in $\curlyO^T$, and therefore by Proposition \ref{prop:endo_degd}, $\curlyO$ is maximally embedded by $\mathbb{Z}[\frac{1+\sqrt{-3}}{2}]$, or $\mathbb{Z}[i]$ respectively. Both of these orders are maximal and of class number $1$, and therefore by the proof of Corollary \ref{cor:ssreduction} $E$ is the reduction modulo $p$ of the unique, up to $\overline{\mathbb{Q}}$-isomorphism, elliptic curve defined over $\overline{\mathbb{Q}}$ with complex multiplication by $\mathbb{Z}[\frac{1+\sqrt{-3}}{2}]$, or $\mathbb{Z}[i]$ respectively, which has $j$-invariant equal to $0$, or $1728$ respectively.
\end{proof}

For the convenience of the reader, we list in Table~\ref{tab:CMorders} below the 13 CM $j$-invariants in $\mathbb{Q}$, associated to the 13 imaginary quadratic orders $O$ of class number $1$ with discriminant $-d=f^2\Delta$, where $\Delta$ is the discriminant of the quadratic field containing $O$.

     \begin{table}[!ht]
        \caption{CM orders \cite[\S 12.C]{Cox} or \cite[Appendix A.3]{SilvermanAdv}}\label{tab:CMorders}
        \centering
   \begin{tabular}{|wc{1em}||wc{4em}|wc{3em}|wc{5em}|wc{3em}|wc{3.5em}|wc{4em}|wc{3em}|}
  \hline
  $j$ & $0$ & $2\cdot 30^3$ & $-3\cdot 160^3$ & $1728$ & $66^3$ & $-15^3$ & $255^3$ \\ \hline
  $O$ & $\mathbb{Z}\left[\frac{1+\sqrt{-3}}{2}\right]$ & $\mathbb{Z}\left[\sqrt{-3}\right]$ & $\mathbb{Z}\left[\frac{3(1+\sqrt{-3})}{2}\right]$ & $\mathbb{Z}\left[\sqrt{-1}\right]$ & $\mathbb{Z}\left[2\sqrt{-1}\right]$ & $\mathbb{Z}\left[\frac{1+\sqrt{-7}}{2}\right]$ & $\mathbb{Z}\left[\sqrt{-7}\right]$ \\ \hline
  $-d$ & $-3$ & $-12$ & $-27$ & $-4$ & $-16$ & $-7$ & $-28$ \\ \hline
  $f$ & $1$ & $2$ & $3$ & $1$ & $2$ & $1$ & $2$ \\ \hline
\end{tabular}

\bigskip

    \begin{tabular}{|wc{1.5em}||wc{3em}|wc{4.5em}|wc{4.5em}|wc{4.5em}|wc{4.5em}|wc{5em}|}
  \hline
  $j$ &  $20^3$ & $-32^3$ & $-96^3$ & $-960^3$ & $-5280^3$ & $-640320^3$ \\ \hline
  $O$ & $\mathbb{Z}\left[\sqrt{-2}\right]$ & $\mathbb{Z}\left[\frac{1+\sqrt{-11}}{2}\right]$ & $\mathbb{Z}\left[\frac{1+\sqrt{-19}}{2}\right]$ & $\mathbb{Z}\left[\frac{1+\sqrt{-43}}{2}\right]$ & $\mathbb{Z}\left[\frac{1+\sqrt{-67}}{2}\right]$ & $\mathbb{Z}\left[\frac{1+\sqrt{-163}}{2}\right]$ \\ \hline
  $-d$ & $-8$ & $-11$ & $-19$ & $-43$ & $-67$ & $-163$ \\ \hline
  $f$ & $1$ & $1$ & $1$ & $1$ & $1$ & $1$ \\ \hline
\end{tabular}
     \end{table}
     
We now end this subsection with a consequence of \Cref{lem:unique_order}, which is a result which appears well known but which we could not find stated in the literature; we include it here for completeness:

\begin{corollary}\label{cor:2isogeny}
   If $E$ is a supersingular elliptic curve defined over $\overline{\mathbb{F}}_p$ with an endomorphism of degree $2$, then $j(E) \in \Fp$. In other words, the vertices with a loop in the supersingular $2$-isogeny graph are curves with $j$-invariants in $\Fp$.
   \end{corollary}
   
\begin{proof}
If $p = 2$ or $3$, there is a unique $\overline{\mathbb{F}}_p$-isomorphism class of supersingular elliptic curve defined over $\overline{\mathbb{F}}_p$ and it has $j$-invariant in $\mathbb{F}_p$ so the statement is vacuously true. (In fact, in both cases the curve does have an endomorphism of degree $2$ since the maximal orders of $B_2$ and $B_3$ each have a unit $u$ of multiplicative order $4$, and the endomorphism corresponding to the element $1+u$ has degree $2$.)

Now let $p \geq 5$, and as usual we write $\End(E) \cong  \curlyO$. Assume that $E$ has an endomorphism of degree $2$ so there exists $\alpha\in \curlyO$ such that $\norm{\alpha}^2=2$ and $\trd(\alpha)\in \mathbb{Z}$. We consider $\beta = 2\alpha - \trd(\alpha)\in \curlyO^T\subseteq \curlyO$ and observe that
\begin{align*}
    0  \leq \norm{\beta}^2 = 4\norm{\alpha}^2-\trd(\alpha)^2 = 8 - \trd(\alpha)^2;
\end{align*}
therefore $\trd(\alpha) = 0, \pm 1, \pm 2$. If $\trd(\alpha)=\pm 2$, then $\norm{\beta}^2=4$ and by Corollary \ref{cor:0and1728}, $j(E) = 1728 \in \mathbb{F}_p$. Otherwise, if $\trd(\alpha)=0,\pm 1$, then $\norm{\beta}^2=7$ or $8$. By Proposition \ref{prop:endo_degd}, $\beta$ corresponds to an embedding of the imaginary quadratic order of discriminant $-\norm{\beta}^2$ into $\curlyO$, which is optimal if and only if $\beta$ is primitive in $\curlyO^T$. 

If $\norm{\beta}^2=7$ and $\beta = n \gamma$ then $\norm{\gamma}^2 = \frac{7}{n^2}$ which is not an integer for any $n \geq 2$, so there can be no such $\gamma$ in $\curlyO^T$ and the embedding is optimal. Since the quadratic order of discriminant $-7$ has class number $1$, by Corollary \ref{cor:ssreduction}, this curve has $j$-invariant in $\mathbb{F}_p$.

If $\norm{\beta}^2=8$ and $\beta = n \gamma$ then $\norm{\gamma}^2 = \frac{8}{n^2}$ which is an integer only if $n=2$ subject to the condition that $n \geq 2$. If $n = 2$ however, $\norm{\gamma}^2 = 2$ which is not congruent to $0$ or $3 \pmod{4}$ so $\gamma$ is not in $\curlyO^T$. Therefore $\beta$ is primitive of norm $8$ and again since the quadratic order of discriminant $-8$ has class number $1$, by Corollary \ref{cor:ssreduction}, this curve has $j$-invariant in $\mathbb{F}_p$.
\end{proof}

\subsection{Rank-$2$ sublattices of determinant $4p$}\label{subsec:rank2sub}

As stated in Proposition \ref{prop:sublat4pkaneko}, in \cite{kaneko1989} Kaneko shows that a necessary condition for a supersingular curve $E$ to have $j$-invariant $j(E)\in \Fp$ is that there exists a rank-$2$ sublattice of $\curlyO^T$ with determinant $4p$. We show that this is also sufficient:

\begin{proposition}\label{prop:sublat4p}
Let $E$ be a supersingular elliptic curve defined over $\overline{\mathbb{F}}_p$. Then $j(E) \in \Fp$ if and only if its Gross lattice $\curlyO^T$ has a rank-$2$ sublattice of determinant $4p$. 
\end{proposition}

\begin{proof}
The proof of the necessity of the condition is given here in Proposition \ref{prop:sublat4pkaneko}.

Conversely, assume that the elements $\beta_1, \beta_2$ form a basis of a rank-$2$ sublattice of $\curlyO^T$ whose determinant is equal to $4p$. Computing the determinant of the Gram matrix directly,
we thus have
\begin{equation*}
\norm{\beta_1}^2\norm{\beta_2}^2-\frac{1}{4}\trd(\beta_1 \overline{\beta}_2)^2=4p.
\end{equation*}
From this, we get $4\norm{\frac{1}{2}\beta_1 \overline{\beta}_2}^2-\trd(\frac{1}{2}\beta_1 \overline{\beta}_2)^2=4p$.

Now set $\alpha= \frac{1}{2}\beta_1 \overline{\beta}_2 - \frac{1}{4} \trd(\beta_1 \overline{\beta}_2 )$. We claim that $\alpha  \in \curlyO$ and that it is of trace $0$ and norm $p$. This implies that $\alpha$ is purely inseparable, and can be factored as $\phi \circ \pi_p$ where $\pi_p \colon E \rightarrow E^{(p)}$ is the $p$th-power Frobenius map and $\phi$ is an isomorphism. It follows that $E$ is $\overline{\mathbb{F}}_p$-isomorphic to its image under $\pi_p$, thus $j(E)= \pi_p(j(E)) = j(E)^p$. Therefore  $j(E) \in \Fp$, and we are done once our claims concerning $\alpha$ are proven. 

To establish these, let $\alpha_1, \alpha_2 \in \curlyO$ be such that $\beta_i= 2 \alpha_i- \trd(\alpha_i)$. Then 
\begin{equation}\label{eq:betas}
\frac{1}{2}\beta_1 \overline{\beta}_2 = 2 \alpha_1 \overline{\alpha}_2 - (\trd(\alpha_2) \alpha_1 + \trd(\alpha_1) \overline{\alpha}_2) + \frac{1}{2} \trd(\alpha_1) \trd(\alpha_2), 
\end{equation}
and 
\begin{equation}\label{eq:tracebetas}
\frac{1}{2}\trd(\beta_1 \overline{\beta}_2)= 2 \trd(\alpha_1 \overline{\alpha}_2) - \trd(\alpha_1) \trd(\alpha_2). 
\end{equation}
Combining equations \eqref{eq:betas} and \eqref{eq:tracebetas} we have 
\begin{align*}
\alpha & = 2\alpha_1\overline{\alpha}_2 - (\trd(\alpha_2) \alpha_1 + \trd(\alpha_1) \overline{\alpha}_2) + \trd(\alpha_1) \trd(\alpha_2) - \trd(\alpha_1 \overline{\alpha}_2).
\end{align*}
This belongs to $\curlyO$ since $\alpha_1, \overline{\alpha}_2$, as well as any integer multiples of them belong to $\curlyO$ and $\trd(\alpha_1),  \trd(\alpha_2) \in \mathbb{Z}$. 

Using the definition of $\alpha$, we also have
\begin{align*}
\trd(\alpha) 
& = \frac{1}{2} \trd(\beta_1\overline{\beta}_2) - \frac{1}{2} \trd(\beta_1\overline{\beta}_2) =0,
\end{align*}
as well as
\begin{align*}
\norm{\alpha}^2 & = \left(\frac{1}{2}\beta_1 \overline{\beta}_2 - \frac{1}{4} \trd(\beta_1 \overline{\beta}_2 )\right)\left(\frac{1}{2}\overline{\beta}_1 \beta_2 - \frac{1}{4} \trd(\beta_1 \overline{\beta}_2 )\right)\\
& = \frac{1}{4}\norm{\beta_1}^2\norm{\beta_2}^2-\frac{1}{16} \trd(\beta_1\overline{\beta}_2)^2 = p.
\end{align*}
\end{proof}

We will in fact need a slighter stronger version of Proposition~\ref{prop:sublat4p}; we present this result here:

\begin{proposition}\label{prop:lalbet12}
Let $E$ be a supersingular elliptic curve defined over $\overline{\mathbb{F}}_p$, and $\{\beta_1, \beta_2, \beta_3\}$ be a successive minimal basis of the Gross lattice $\curlyO^T$ of $E$. Then $j(E) \in \Fp$ if and only if the sublattice $\Lambda=\langle \beta_1, \beta_2 \rangle$ of $\curlyO^T$ has determinant $4p$.
 \end{proposition}
 
 \begin{proof}
  To prove this result, by Proposition~\ref{prop:sublat4p}, it suffices to show that if $j(E)\in \Fp$ then $\det(\Lambda)=4p$. Hence, suppose that $j(E)\in \Fp$  and assume for a contradiction that $\det(\Lambda)\ne 4p$.  By Proposition \ref{prop:sublattice4np}, $\det(\Lambda)$ is a positive integer multiple of $4p$, and thus we must have $\det(\Lambda)\ge 8p$. 

   By Proposition~\ref{prop:sublat4p}, there is a sublattice $\Lambda'$ of $\curlyO^T$ such that $\det(\Lambda')=4p$; let $\{\gamma_1, \gamma_2\}$ be a successive minimal basis for $\Lambda'$. Then we have that $\Lambda'=\langle\gamma_1, \gamma_2\rangle$, $\|\beta_1\|^2 \le \norm{\gamma_1}^2$ and $\|\beta_2\|^2 \le \|\gamma_2\|^2$. Moreover, by equation \eqref{eq:hermite}
\begin{equation}\label{eq:boundgamma}
    \|\gamma_1\|^4 \le \|\gamma_1\|^2 \|\gamma_2\|^2 \le \frac{4}{3} \det(\Lambda')= \frac{16}{3} p.
\end{equation}

Denote by $\{\beta_1, \beta_2^*\}$  and $\{\gamma_1, \gamma_2^*\}$, respectively, the Gram-Schmidt orthogonalization of the bases $\{\beta_1, \beta_2\}$ and $\{\gamma_1, \gamma_2\}$, respectively. Then $\det(\Lambda)= \norm{\beta_1}^2 \norm{\beta^*_2}^2$, and hence 
\begin{equation}\label{eq:boundbeta2}
    \|\beta_2\|^2 \ge \|\beta_2^*\|^2= \frac{\det(\Lambda)}{\|\beta_1\|^2} \ge \frac{8p}{\|\gamma_1\|^2}.
\end{equation}

Now let $\mu_{2,1}= \frac{(\gamma_2, \gamma_1)}{\norm{\gamma_1}^2}$, we have $|\mu_{2,1}| \le \frac{1}{2}$ by \Cref{lem:size_reduced}. Since $4p= \det(\Lambda')= \|\gamma_2^*\|^2 \|\gamma_1\|^2$, we obtain that 
\begin{equation}\label{eq:boundbeta2_2}
    \|\beta_2\|^2 \le \|\gamma_2\|^2= \|\gamma_2^*\|^2+(\mu_{2,1})^2 \|\gamma_1\|^2 \le \frac{4p}{\|\gamma_1\|^2} +\frac{1}{4} \|\gamma_1\|^2.
\end{equation}

     From inequalities \eqref{eq:boundbeta2} and \eqref{eq:boundbeta2_2}, we have  
     \[\frac{8p}{\|\gamma_1\|^2} \le \frac{4p}{\|\gamma_1\|^2} +\frac{1}{4} \|\gamma_1\|^2.\]
     This implies that $\|\gamma_1\|^4 \ge 16p$ which is a contradiction by equation \eqref{eq:boundgamma}.
 \end{proof}

\begin{remark}
    If $j(E) \in \mathbb{F}_p$, it is not the case that every rank-2 sublattice of determinant $4p$ in $\mathcal{O}^T$ must have a basis attaining the first two successive minima of $\curlyO^T$. Indeed if $p \neq 3$, and $\{\beta_1,\beta_2,\beta_3\}$ is a successive minimal basis for the Gross lattice of a supersingular elliptic curve with $j$-invariant equal to $1728$, we have that $\det(\langle \beta_1,\beta_3\rangle) = 4p$, as demonstrated in Proposition \ref{prop:1728} below, and this sublattice does not contain any element of norm $p$ since if $\beta = x_1 \beta_1+ x_3 \beta_3$, then $\norm{\beta}^2 = x_3^2p + (x_3+2x_1)^2$.\end{remark}

    \begin{remark}\label{rem:CGconverse}
    As we recalled in Lemma \ref{lem:CGd1d2}, \cite{Chevyrev-Galbraith} shows that if $j(E) \in \mathbb{F}_p$, then $D_1D_2 < \frac{16p}{3}$. If conversely $D_1D_2 < \frac{16p}{3}$, using equation \eqref{eq:hermite}, we have then that $\det(\langle \beta_1, \beta_2\rangle) < \frac{16p}{3} < 8p$, where as usual $\beta_1$ and $\beta_2$ are the first two vectors of a successive minimal basis of $\curlyO^T$. By Proposition \ref{prop:sublattice4np}, this forces $\det(\langle \beta_1, \beta_2\rangle) = 4p$ and by Proposition~\ref{prop:lalbet12} we then have that $j(E) \in \mathbb{F}_p$. Hence we conclude that in fact $D_1D_2 < \frac{16p}{3}$ if and only if $j(E) \in \mathbb{F}_p$.
\end{remark}

\subsection{Successive minimal bases for $j$-invariants $0$ and $1728$}\label{subsec:0and1728}
In this section, we compute a normalized successive minimal basis, and its Gram matrix, for the Gross lattices of elliptic curves with $j$-invariants $0$ and $1728$ when they are supersingular. 

\begin{proposition}\label{prop:0}
Let $p\equiv 2 \pmod 3$ be a prime, $\mi,\mj,\mk$ be elements of $B_p$ such that $\mi^2=-3, \mj^2=-p, \mk= \mi\mj$, and $E_0$ be the elliptic curve defined over $\overline{\mathbb{F}}_p$ with $j(E_0)=0$. Then $E_0$ is supersingular and the Gross lattice $\curlyO^T$ of $E_0$ has a normalized successive minimal basis given by 
\begin{equation*}
\left\{\mi, \frac{\mi+3\mj-\mk}{3}, \frac{-\mi-2\mk}{3} \right\}=: \{\beta_1, \beta_2, \beta_3\}.
\end{equation*}
Moreover, the Gram matrix of this basis is
\[\begin{pmatrix}
3 & 1 & 1 \\
1 & \frac{4p+1}{3} & -\frac{2p-1}{3} \\
1 & - \frac{2p-1}{3} & \frac{4p+1}{3}
\end{pmatrix}.\]
\end{proposition}

\begin{proof}
Let $\curlyO$ be a maximal order in the quaternion algebra $B_p$ such that $ \curlyO \cong \End(E_0)$.  Denote by $\{\beta_1, \beta_2, \beta_3\}$ a normalized successive minimal basis of $\curlyO^T$ and as usual set $D_i= \|\beta_i\|^2$. 
    
In this case, $\curlyO$ is maximally embedded by $\mathbb{Z}[\frac{1+\sqrt{-3}}{2}]$, see for example \cite[Theorem 10.7]{washington}. For $\alpha \in \curlyO$ corresponding to $\frac{1+\sqrt{-3}}{2}$, we have $\trd(\alpha)=1$ and $\norm{\alpha}^2=1$, and setting $\beta_1= 2\alpha-1 \in \curlyO^T$, we have $\norm{\beta_1}^2=3$ and $D_1=3$ is the first successive minimum of $\curlyO^T$ since elements of the Gross lattice have norm congruent to $0$ or $3$ modulo $4$. 

The Gram matrix of the basis $\{\beta_1, \beta_2, \beta_3\}$ of $\curlyO^T$ is thus
 \[G_{\{\beta_1,\beta_2,\beta_3\}}=\begin{pmatrix}
3 & t_{12} & t_{13} \\
t_{12} & D_2 & t_{23} \\
t_{13} & t_{23} & D_3
\end{pmatrix},
\]
where $t_{12}, t_{13}, t_{23} \in \mathbb{Z}$, and $0 \leq t_{12},t_{13}  \le \frac{D_1}{2}=\frac{3}{2}$ by Lemma~\ref{lem:Gram_matrix}. 
By Proposition~\ref{prop:lalbet12}, the lattice $\Lambda_1=\langle\beta_1, \beta_2\rangle$ has determinant $4p$, or $4p= 3D_2 -t_{12}^2$.
This forces $t_{12}=1$ since $p\equiv 2 \pmod 3$, and hence $D_2= \frac{4p+1}{3}$.

It is known that 
$\curlyO = \langle 1, \frac{1+\mi}{2}, \frac{\mj-\mk}{2}, \frac{\mi-\mk}{3} \rangle$ for $\mi^2=-3$ and $\mj^2=-p$ by \cite{Ibukiyama}. We can compute $\curlyO^T$ using this presentation and obtain the basis $\{\mi , \mj-\mk, \frac{2(\mi-\mk)}{3} \} =: \{ \gamma_1,\gamma_2,\gamma_3\}$ for 
$\curlyO^T$. We obtain a smaller basis for $\curlyO^T$ by performing unimodular $\mathbb{Z}$-operations as follows:
First, change the sign of $\gamma_3$, then add $\gamma_1$ to $\gamma_3$ to obtain the new basis
$\{\mi, \mj-\mk, \frac{\mi+2\mk}{3}\}=: \{\delta_1, \delta_2, \delta_3\}$. For this new basis, add $\delta_3$ to $\delta_2$ which yields the basis 
$\{ \mi, \frac{\mi+3\mj-\mk}{3}, \frac{\mi+2\mk}{3} \}=: \{\beta_1, \beta_2, \beta_3\}$. Note that $\|\beta_1\|^2=3=D_1$ and the last two elements have the same norm, i.e., $\|\beta_2\|^2= \|\beta_3\|^2 = \frac{4p+1}{3}= D_2$, which implies $D_3= \frac{4p+1}{3}$, and the basis $\{\beta_1, \beta_2, \beta_3\}$ attains the successive minima of $\curlyO^T$. The result is then obtained by computing the Gram matrix of this basis. 
\end{proof}

\begin{proposition}\label{prop:1728}
Let $p\equiv 3 \pmod 4$ be a prime, $\mi,\mj,\mk$ be elements of $B_p$ such that $\mi^2=-1, \mj^2=-p, \mk= \mi\mj$, and $E_{1728}$ be the elliptic curve defined over $\overline{\mathbb{F}}_p$ with $j(E_{1728})=1728$. Then $E_{1728}$ is supersingular and if $p>3$, then the Gross lattice $\curlyO^T$ of $E_{1728}$ has a normalized successive minimal basis given by
\begin{equation*}
\{ 2\mi,\mj, \mi-\mk \}=: \{\beta_1, \beta_2, \beta_3\}.
\end{equation*}
Moreover, the Gram matrix of this basis is
\[ 
\begin{pmatrix}
4 & 0 & 2 \\
0 & p & 0 \\
2 & 0 & p+1
\end{pmatrix}.\]
\end{proposition}
\begin{proof}
Let $\curlyO$ be a maximal order in the quaternion algebra $B_p$ such that $ \curlyO \cong \End(E_{1728})$.  Denote by $\{\beta_1, \beta_2, \beta_3\}$ a normalized successive minimal basis of $\curlyO^T$, and again set $D_i= \|\beta_i\|^2$. 
    
This time $E_{1728}$ is maximally embedded by $\mathbb{Z}[\sqrt{-1}]$, see for example \cite[Theorem 10.7]{washington}. There is therefore $\alpha \in \curlyO$ with $\trd(\alpha)=0$ and $\norm{\alpha}^2=1$, and setting $\beta_1= 2\alpha \in \curlyO^T$, we have $\norm{\beta_1}^2 = 4$. We now argue that no element $\gamma$ of $\curlyO^T$ has $\norm{\gamma}^2 = 3$ to conclude that $D_1 = 4$ in this case. Indeed, by Corollary \ref{cor:0and1728}, if there were two elements of respective norms $3$ and $4$ in $\curlyO^T$, then we would have $0 \equiv 1728 \pmod{p}$, which implies $p =2$ or $p=3$, and these primes are excluded from consideration in this proposition. Therefore $D_1 = 4$.

The Gram matrix of the basis $\{\beta_1, \beta_2, \beta_3\}$ of $\curlyO^T$ is thus
\[G_{\{\beta_1,\beta_2,\beta_3\}}=\begin{pmatrix}
4 & t_{12} & t_{13} \\
t_{12} & D_2 & t_{23} \\
t_{13} & t_{23} & D_3
\end{pmatrix},
\]
where $t_{12},t_{13},t_{23} \in \mathbb{Z}$ and $ 0\le t_{12}, t_{13}  \le \frac{D_1}{2}=2$ by Lemma~\ref{lem:Gram_matrix}. We again use the fact that the lattice $\Lambda_1=\langle\beta_1, \beta_2\rangle$ has determinant $4p$ by Proposition~\ref{prop:lalbet12}, so $4p= 4D_2 -t_{12}^2$, and
\begin{equation}\label{eq:D2_1728}
    D_2= \frac{4p+t_{12}^2}{4} \quad \text{for} \quad  t_{12}\in \{0, 2\}.
\end{equation}

As shown in \cite{Ibukiyama}, $\curlyO = \langle 1, \mi, \frac{1+\mj}{2},  \frac{\mi+\mk}{2} \rangle$  for $\mi^2=-1$ and $\mj^2=-p$. Using this presentation to compute $\curlyO^T$ we obtain the basis $\{2\mi, \mj, \mi+\mk \} =:\{ \gamma_1,\gamma_2,\gamma_3\}$ for $\curlyO^T$. Replacing $\gamma_3$ with $\gamma_1-\gamma_3$ gives the basis $\{2\mi, \mj, \mi-\mk \} =:\{ \delta_1,\delta_2,\delta_3\}$. 
Here $\|2\mi\|^2=4=D_1$, and the last two elements in the basis have norm $\|\mj\|^2= p$ and $\|\mi-\mk\|^2=p+1$, so we must have $D_2\le p$ and $D_ 3\le p+1$. Using equation \eqref{eq:D2_1728}, we have $D_2=p$, $t_{12}=0$, and $p\le D_3\le p+1$.

Using the above $\mathbb{Z}$-basis of $\curlyO^T$, a straightforward computation shows that the only elements of norm $p$ in $\curlyO^T$ are $\pm \mj$. As these two vectors are linearly dependent, we conclude that $D_3> p$, and therefore $D_3= p+1$, and the basis $\{2\mi, \mj, \mi-\mk\}$ attains the successive minima of $\curlyO^T$.

Now, the lattice $\Lambda_2=\langle\beta_1, \beta_3\rangle$ has determinant $4np$ for some positive integer $n$ by Proposition \ref{prop:sublattice4np}. Thus, $4np= 4(p+1) - t_{13}^2$, which forces $n=1$ and $t_{13}=2$.
Finally, since $4p^2=\det(\curlyO^T)= 4p^2-4t_{23}^2$, we have that $t_{23}=0$, and the result follows.
\end{proof}

\begin{remark}\label{rem:1728}
Since our result does not cover the case of $p=3$, for completeness we work it out here. If $p = 3$ then $1728 \equiv 0 \pmod{p}$, so $\curlyO$ is maximally embedded by $\mathbb{Z}[\frac{1+\sqrt{-3}}{2}]$ by \cite[Theorem 10.7]{washington} and $D_1 = 3$. By the proof of Proposition \ref{prop:1728}, it follows then that $3 \leq D_2 \leq 4$ since $\curlyO^T$ contains an element of norm $4$, which is necessarily linearly independent from the element of norm $3$. Now fix $\{\beta_1,\beta_2,\beta_3\}$ any normalized successive minimal basis of $\mathcal{O}^T$ and write $t_{ij} = (\beta_i,\beta_j)$. Since $D_2 = \frac{t_{12}^2+12}{3}$ with $0 \leq t_{12} \leq 1$, we have $t_{12}=0$ and $D_2 = 4$. Noting that $3 D_3 -t_{13}^2$ is a multiple of $12$ and $t_{13} \in \{0,1\}$, we conclude that $t_{13} = 0$. Finally, using Lemma \ref{lem:detOT}, $\det(\curlyO^T) =  36$, so $4D_3-t_{23}^2 = 12$ with $|t_{23}| \leq 2$ and $D_3 \geq 4$, so $D_3 = 4$ and $t_{23}=\pm 2$. Hence the Gram matrix of a normalized successive minimal basis $\{\beta_1,\beta_2,\beta_3\}$ in this case is given by
\[ G_{\{\beta_1,\beta_2,\beta_3\}} = \begin{pmatrix}
3 & 0 & 0 \\
0 & 4 & \pm 2 \\
0 & \pm 2 & 4
\end{pmatrix}. 
\]
\end{remark}

\subsection{The third successive minimum of Gross lattices when $j(E)\notin \Fp$}\label{subsec:D3notinFp}
   
In this section, we prove the first part of our main result that characterizes the field of definition of the $j$-invariant of a supersingular elliptic curve by the value of the third successive minimum of its Gross lattice. We begin by establishing that when $j(E) \in \mathbb{F}_{p^2} \setminus \mathbb{F}_p$, we have $D_3 < p$, and then refine this bound, 
determining the extent to which $D_3$ must be less than $p$ when $j(E) \in \mathbb{F}_{p^2} \setminus \mathbb{F}_p$.

\begin{lemma}\label{lem:EnotinFp_D3lessthanp}
    Let $E$ be a supersingular curve defined over $\overline{\mathbb{F}}_p$ and $D_3$ be the third successive minimum of its Gross lattice. Then $j(E) \in \mathbb{F}_p$ if and only if $D_3 \geq p$.
\end{lemma}

\begin{proof}
First suppose that $j(E) \in \mathbb{F}_p$, and let $\{\beta_1,\beta_2,\beta_3\}$ be a successive minimal basis for the Gross lattice of $E$ and $\{\beta_1,\beta^*_2,\beta^*_3\}$ be its Gram-Schmidt orthogonalization. By Proposition~\ref{prop:lalbet12} we have
\begin{equation*}
\norm{\beta_1}^2 \norm{\beta^*_2}^2 = \det(\langle \beta_1, \beta_2 \rangle) = 4p,
\end{equation*}
and by Lemma \ref{lem:detOT}
\begin{equation*}
\norm{\beta_1}^2 \norm{\beta^*_2}^2\norm{\beta_3^*}^2 = \det(\curlyO^T) = 4p^2,
\end{equation*}
and therefore
\begin{equation*}
D_3 \geq \norm{\beta_3^*}^2 = p.
\end{equation*}

Now suppose that $j(E) \in \mathbb{F}_{p^2} \setminus \mathbb{F}_p$, and again let $\{\beta_1, \beta_2, \beta_3\}$ be a normalized successive minimal basis of $\curlyO^T$ with $\|\beta_i\|^2= D_i$ for $i = 1,2,3$. We recall first from equation \eqref{eq:prodbound} that 
\begin{equation}\label{eq:d1d2d3}
4p^2 \leq D_1 D_2 D_3 \leq 8p^2.
\end{equation}
Furthermore, by Proposition \ref{prop:sublattice4np} and Proposition~\ref{prop:lalbet12}, since $j(E)\notin \Fp$, we have that $\det ( \langle \beta_1,\beta_2 \rangle) = 4np$ for $n \geq 2$ an integer. Using the Hermite bound (equation \eqref{eq:hermite}) for this sublattice $\langle \beta_1,\beta_2 \rangle$ we have
\begin{equation}\label{eq:d1d2}
4np \leq D_1 D_2 \leq \frac{16}{3}np.
\end{equation}
Combining equations \eqref{eq:d1d2d3} and \eqref{eq:d1d2}, we have 
\begin{equation*}
D_3 \leq \frac{8p^2}{4np} \leq p
\end{equation*}
since $n \geq 2$. 

We now show that $D_3 \neq p$: Indeed, if $D_3 = p$, we must have $n = 2$ 
and again combining equations \eqref{eq:d1d2d3} and \eqref{eq:d1d2}, $D_1 D_2 = 8p$. If $p \leq 31$, then every supersingular elliptic curve defined over $\overline{\mathbb{F}}_p$ has $j$-invariant contained in $\mathbb{F}_p$, hence we may assume that $p \geq 37$. Since $D_1,D_2 \le D_3=p$ and $p \geq 37$, it follows that $D_1= 8$ and $D_2=p$. However, by Corollary \ref{cor:ssreduction}, since the imaginary quadratic order of discriminant $-8$ has class number $1$, this implies that the $j$-invariant of the curve under consideration is in $\mathbb{F}_p$. Hence if $j(E) \in \mathbb{F}_{p^2} \setminus \mathbb{F}_p$, we have $D_3<p$.
\end{proof}

To refine our bound on $D_3$ further, we will need the following lemmata.

\begin{lemma}\label{lem:upperboundD1}
Let $E$ be a supersingular curve defined over $\overline{\mathbb{F}}_p$ and $\{\beta_1,\beta_2,\beta_3\}$ be a successive minimal basis of $\curlyO^T$ with $\det(\langle \beta_1, \beta_2 \rangle)= 4np$ for $n$ a positive integer. Then if $Y>0$ is a real number and $p > \frac{16n}{3Y^2}$, we have $D_1=\norm{\beta_1}^2< Y p$. 
\end{lemma}
\begin{proof}
Assume by contradiction that $D_1 \ge  Y p$. Applying the Hermite bound given in equation \eqref{eq:hermite} to the sublattice $\langle \beta_1, \beta_2\rangle$, we have 
\begin{equation*}
4np \leq D_1 D_2 \leq \frac{16}{3}np,
\end{equation*}
and since $D_1 \le D_2$, we obtain that $D_1^2 \le \frac{16np}{3}$. Thus, $Y^2 p^2 \le \frac{16np}{3}$ and $p \le \frac{16n}{3Y^2}$, contradicting the assumption. 
\end{proof}

The following lemma is a straightforward calculus problem:
\begin{lemma}\label{lem:max}
    Let $a, b$ be positive real numbers such that $b \ge \sqrt{a}$, and let $f(x)= x +\frac{a}{x}$. Then $f(x)$ is increasing on the interval $[\sqrt{a}, b]$, and in particular, 
    \begin{equation*}
    \max_{[\sqrt{a}, b]} f(x) = f(b) = b + \frac{a}{b}.
    \end{equation*}
\end{lemma}

We can now prove our upper bound on the value of the third successive minimum of the Gross lattice when the supersingular elliptic curve has $j$-invariant in $\Fpp\backslash \Fp$.

\begin{theorem}\label{thm:D3notoverFp}
Let $E$ be a supersingular curve defined over $\overline{\mathbb{F}}_p$ with $j(E) \in \Fpp\backslash \Fp$, and $D_3$ be the third successive minimum of its Gross lattice. Then 
\begin{equation*}
D_3 \le \frac{3}{5}p +5.
\end{equation*} 
\end{theorem}

\begin{proof}
Throughout, we use the notation developed in Lemma~\ref{lem:EnotinFp_D3lessthanp} and its proof. Since $j(E)\notin \Fp$, we have that $\det ( \langle \beta_1,\beta_2 \rangle) = 4np$ for $n \geq 2$ an integer by Proposition \ref{prop:sublattice4np} and Proposition~\ref{prop:sublat4p}. If $n \geq 4$, this immediately yields that $D_3 \leq \frac{8p^2}{16p} = \frac{p}{2} < \frac{3}{5}p+5$ by equations \eqref{eq:d1d2d3} and \eqref{eq:d1d2}, and we are done.

Carrying on in the case of $n = 2$ or $3$, we first consider the case of $p$ large, and let $p\ge 67$. By Lemma \ref{lem:size_reduced}, if $\{ \beta_1, \beta_2, \beta_3\}$ is a successive minimal basis for $\curlyO^T$, then the pairs $\{ \beta_1, \beta_2\}$ and $\{\beta_1,\beta_3\}$ are size-reduced, but the pair $\{\beta_2,\beta_3\}$ may not be. We therefore proceed to size-reduce the basis $\{ \beta_1, \beta_2, \beta_3\}$ to obtain the new basis $\{ \beta_1, \beta_2, \beta'_3\}$ for $\curlyO^T$, where, explicitly,
 \begin{equation*}
 \beta'_{3} = \beta_3 - \round{\mu_{3,2}} \beta_2 - \round{\mu_{3,1} - \round{\mu_{3,2}} \mu_{2,1}} \beta_1.
 \end{equation*}
Furthermore, we write
\begin{equation*}
\mu'_{3,1} = \frac{(\beta'_3,\beta_1)}{\norm{\beta_1}^2} \quad \text{and} \quad \mu'_{3,2} = \frac{(\beta'_3,\beta^*_2)}{\norm{\beta^*_2}^2}.
\end{equation*}

By the properties of size-reduced bases, we have $|\mu'_{3,1}| \le \frac{1}{2}$ and  $|\mu'_{3,2}| \le \frac{1}{2}$ as well as  
\begin{equation}\label{eq:beta3p0}
\norm{\beta'_3}^2 = \norm{\beta^*_3}^2 + (\mu'_{3,1})^2 \norm{\beta_1}^2 + (\mu'_{3,2})^2 \norm{\beta_2^*}^2,
\end{equation}
where $\{\beta_1,\beta_2^*, \beta_3^*\}$ is the common Gram-Schmidt orthogonalization of the two bases $\{ \beta_1, \beta_2, \beta_3\}$ and $\{ \beta_1, \beta_2, \beta'_3\}$.

Since $\beta_3'$ is linearly independent from $\beta_1$ and $\beta_2$ and using equation \eqref{eq:beta3p0}, we therefore have
\begin{align}\label{eq:d3bound1}
D_3 \leq \norm{\beta'_3}^2 \leq \norm{\beta^*_3}^2  + \frac{1}{4}\left(\norm{\beta_1}^2 + \norm{\beta_2^*}^2\right). 
\end{align}

Using equation \eqref{eq:d3bound1} and the relation $\norm{\beta_1}^2 \norm{\beta^*_2}^2 = \det ( \langle \beta_1,\beta_2 \rangle) =4np$, we can further write this as
\begin{equation}\label{eq:d3bound}
D_3 \leq \norm{\beta^*_3}^2  + \frac{1}{4}\left(X + \frac{4np}{X}\right),
\end{equation}
where $X=\max\{\norm{\beta^*_2}^2, \norm{\beta_1}^2 \} \ge \sqrt{4np}$.

We now turn our attention specifically to the case of $n = 3.$ In this case 
\begin{equation*}
\norm{\beta^*_3}^2 = \frac{\det(\curlyO^T)}{\det ( \langle \beta_1,\beta_2 \rangle)}=\frac{4p^2}{12p} = \frac{1}{3}p,
\end{equation*} 
and we note that 
\begin{equation*}
X=\max\{\norm{\beta^*_2}^2, \norm{\beta_1}^2 \}  \leq \norm{\beta_2}^2  \leq \norm{\beta_3}^2 = D_3 < p
\end{equation*}
by Lemma~\ref{lem:EnotinFp_D3lessthanp}. Since $p \geq 13$, we may apply Lemma~\ref{lem:max} with $a = 12p$, and $b= p$, and so  $X + \frac{12p}{X} \leq p + \frac{12p}{p}= p+12$. Equation~\eqref{eq:d3bound} then gives
\begin{equation*}
D_3 \leq \frac{p}{3} + \frac{1}{4}\left(p+12\right) = \frac{7}{12}p + 3 < \frac{3}{5}p+5,
\end{equation*}
and we are done.

Finally, we consider the case of $n = 2$, and now $\norm{\beta^*_3}^2 = \frac{4p^2}{8p} = \frac{1}{2}p$. 
We further consider two subcases: 

First, suppose that $X \leq \frac{2}{5}p$. In this case, since $p \geq 53$ we may apply Lemma~\ref{lem:max} with $a = 8p,$ and $b= \frac{2}{5}p$, and obtain that $X + \frac{8p}{X} \leq \frac{2p}{5} + \frac{8p}{2p/5}= \frac{2p}{5} + 20$. In this case, equation~\eqref{eq:d3bound} gives 
\begin{equation*}
D_3 \leq \frac{p}{2} + \frac{1}{4}\left(\frac{2p}{5}+20\right) = \frac{3}{5}p + 5.
\end{equation*}

Therefore, it remains to consider the last case, for which we recall that we have $\det (\langle \beta_1,\beta_2\rangle) = 8p$ and $X = \max\{\norm{\beta^*_2}^2, \norm{\beta_1}^2 \} > \frac{2}{5}p$. By \Cref{lem:upperboundD1}, we have that for $p\ge 67$, $\norm{\beta_1}^2 < \frac{2}{5}p,$ and therefore $\norm{\beta^*_2}^2=X > \frac{2}{5}p$. 
In this case, $\det (\langle \beta_1,\beta_2\rangle)$ is small and $\beta^*_2$ is long, which forces $\beta_1$ to be very short. Indeed, since $\beta_1$ and $\beta^*_2$ are orthogonal, we have
\begin{equation}
8p = \det(\langle \beta_1, \beta_2 \rangle)= \norm{\beta_1}^2 \norm{\beta^*_2}^2 > \frac{2}{5}p \norm{\beta_1}^2,
\end{equation}
or $D_1 = \norm{\beta_1}^2 < 20$. Moreover, by Lemma~\ref{lem:EnotinFp_D3lessthanp}, $\norm{\beta^*_2}^2 \le D_2 \le D_3 <p$, thus $\norm{\beta_1}^2> 8$. 
Since $D_1 \equiv 0, 3 \pmod 4$, we must have $D_1 \in \{11, 12, 15, 16, 19\}$. 

Since $D_1$ is the norm of a shortest nonzero element $\beta_1$ of $\curlyO^T$, $\beta_1$ is a primitive vector of $\curlyO^T$. Therefore by Proposition \ref{prop:endo_degd}, $\curlyO$ is maximally embedded by the imaginary quadratic order of discriminant $d=-D_1$. If  $d \in \{-11, -12, -16, -19\}$, the quadratic order of discriminant $d$ has class number $1$, and by Corollary \ref{cor:ssreduction}, this elliptic curve has $j$-invariant in $\mathbb{F}_p$.
    
Thus, the only case left to consider for large primes has $D_1=15$, $\det (\langle \beta_1,\beta_2\rangle) = 8p$, and $\norm{\beta^*_3}^2 = \frac{p}{2}$. To handle it, we consider the sublattice $\Lambda= \langle \beta_1, \beta_3'\rangle$ of $\curlyO^T$, with $\det(\Lambda)= 4sp$ for some integer $s\ge 2$ by Proposition \ref{prop:sublattice4np} and Proposition~\ref{prop:sublat4p}. Let $\{\beta_1, \beta^{**}_3\}$ be the Gram-Schmidt orthogonalization of $\{\beta_1, \beta'_3\}$ with Gram-Schmidt coefficient $\mu'_{3, 1}$. Then $\beta^{**}_3 = \beta'_3- \mu'_{3,1} \beta_1$ and 
\begin{equation}\label{eq:beta3p1}
\norm{\beta'_3}^2= \norm{\beta^{**}_3}^2 + (\mu'_{3,1})^2 \norm{\beta_1}^2.
\end{equation}
In addition, we have
\begin{equation}\label{eq:beta3pp}
    \norm{\beta^{**}_3}^2= \frac{\det(\Lambda)}{\norm{\beta_1}^2} = \frac{\det(\Lambda)}{15}.
\end{equation}
Subtracting equation \eqref{eq:beta3p0} from equation \eqref{eq:beta3p1} yields
$$0= \norm{\beta^{**}_3}^2-\norm{\beta^*_3}^2 - (\mu'_{3,2})^2 \norm{\beta^*_2}^2.$$
The last equality, together with equation \eqref{eq:beta3pp}, imply that 
\begin{equation}\label{eq:mu3pp1}
    (\mu'_{3,2})^2 = \frac{\norm{\beta^{**}_3}^2-\norm{\beta^*_3}^2}{\norm{\beta^*_2}^2}=\frac{2\det(\Lambda)- 15 p}{30\norm{\beta^*_2}^2}.
\end{equation}
Since  $|\mu'_{3,2}| \le \frac{1}{2}$ and $15 \norm{\beta_2^*}^2=\norm{\beta_1}^2 \norm{\beta_2^*}^2= 8p$, from equation \eqref{eq:mu3pp1} we obtain 
\begin{equation*}\label{eq:mu3pp2}
    \frac{1}{4} \ge (\mu'_{3,2})^2 = \frac{2\det(\Lambda)- 15 p}{16p},
\end{equation*}
which forces $\det(\Lambda) =8p$ since $\det(\Lambda)$ is a multiple of $4p$ that is strictly greater than $4p$. By equation \eqref{eq:beta3pp}, we must then have $\norm{\beta^{**}_3}^2= \frac{\det(\Lambda)}{\norm{\beta_1}^2} = \frac{8p}{15}$, and applying equation \eqref{eq:beta3p1}, we then have 
   $$D_3 \leq \norm{\beta'_3}^2\le \frac{8p}{15} + \frac{15}{4}\le \frac{3p}{5} + 5,$$
which completes the proof in this case. 

Finally, we handle the small values of $p$ excluded from previous consideration, i.e., $ p<67$. If $p \leq 31$, every supersingular curve defined over $\overline{\mathbb{F}}_p$ has $j$-invariant in $\mathbb{F}_p$, and there is nothing to prove. For each $37 \leq p<67$, we compute a $\mathbb{Z}$-basis for each maximal order $\curlyO$ in $B_p$, and then a successive minimal basis for the Gross lattice $\curlyO^T$.\footnote{\url{https://github.com/gkorpal/minimal-gross/blob/main/scripts/finite_cases/ref/cases_1_100.txt}} By Lemma~\ref{lem:EnotinFp_D3lessthanp}, the ones with $D_3<p$ correspond to curves with $j(E) \in \Fpp\backslash \Fp$, and we can verify that all of these curves have $D_3 \le \frac{3p}{5}+5$, which completes the proof.
\end{proof}

The following proposition shows that the upper bound $D_3\le \frac{3}{5}p +5$ of Theorem~\ref{thm:D3notoverFp} is quite tight. Indeed, when $p\equiv 13, 17 \pmod{20}$ and $p \ge 113$, there exists a supersingular elliptic curve with $j$-invariant in $\Fpp\setminus\Fp$ whose Gross lattice has third successive minimum $\frac{3}{5}p+\frac{1}{5} \leq D_3$:

\begin{proposition}\label{prop:Fp2_D1is20}
     Let $p \ge 113$ be a prime with $p\equiv 13,17\pmod {20}$. There exists a supersingular elliptic curve with $j$-invariant in $\Fpp\setminus\Fp$, such that the Gram matrix of any normalized successive minimal basis $\{\beta_1,\beta_2,\beta_3\}$ of its Gross lattice $\curlyO^T$ is
    \[G_{\{\beta_1,\beta_2,\beta_3\}} = \begin{pmatrix}
        20 & 2r & 2s \\
        2r & \frac{2p+r^2}{5} & \frac{- p + rs}{5} \\
        2s & \frac{- p + rs}{5} & \frac{3p+s^2}{5}
    \end{pmatrix}\]
    where $(r,s) =(3,1)$ if $p\equiv 13 \pmod {20}$, and $(r,s)=(1, 2)$ if $p\equiv 17 \pmod {20}$.       
    \end{proposition}
\begin{proof}
Let $E$ be an elliptic curve defined over $\overline{\mathbb{Q}}$ and which has complex multiplication by the imaginary quadratic order of discriminant $-20$.\footnote{There are in fact two such curves up to $\overline{\mathbb{Q}}$-isomorphism; an example of one of them has LMFDB elliptic curve label \href{https://www.lmfdb.org/EllipticCurve/2.2.5.1/4096.1/n/1}{2.2.5.1-4096.1-n1}.} Such an elliptic curve has $j$-invariant in $\mathbb{Q}(\sqrt{5})$, and a model defined over $\mathbb{Q}(\sqrt{5})$ with good reduction at every prime above $p \neq 2$. Therefore, if $p\equiv 13,17\pmod {20}$, then $p$ is inert in $\mathbb{Q}(\sqrt{5})$ and the reduction of $E$ modulo the prime ideal $(p)$ above $p$ in $\mathbb{Q}(\sqrt{5})$ is a supersingular elliptic curve $\widetilde{E}$ with a model defined over $\Fpp$, the quotient of the ring of integers of $\mathbb{Q}(\sqrt{5})$ by the prime ideal $(p)$. By Proposition \ref{prop:endo_degd}, the Gross lattice $\curlyO^T$ of $\widetilde{E}$ has an element of norm $20$, and by \cite[Lemma 2.2.1]{Goren-Lauter}, by choosing $p\ge 113$, we can ensure that the first successive minimum of $\curlyO^T$ is $D_1=20$, or in other words that the endomorphism ring of $\widetilde{E}$ is not maximally embedded by any other imaginary quadratic order of larger discriminant.

We first show that $j(\widetilde{E})\notin \Fp$. By way of contradiction, let $\{\beta_1, \beta_2, \beta_3\}$ be a normalized successive minimal basis of $\curlyO^T$, and suppose that $j(\widetilde{E})\in \Fp$, then $\det\langle\beta_1, \beta_2 \rangle = 4p$ by Proposition~\ref{prop:lalbet12}. Using the same notation as in the proof of Proposition~\ref{prop:0}, by \Cref{lem:size_reduced} there is an integer $x$ with $0\le x \le 10$ such that  $20 D_2-x^2=4p$. Thus, $x= 2r$ for some integer $r$ with $0\le r \le 5$ and $5|(p+r^2)$. If $p\equiv 13 \pmod{20},$ we then have $r^2\equiv 2 \pmod 5$ which is impossible. The case of $p\equiv 17 \pmod{20}$ is handled similarly. Thus, we must have that $j(\widetilde{E})\notin \Fp$.

We now compute the Gram matrix of the normalized successive minimal basis. With the same notation as above, we now have that $20 D_2-x^2=\det\langle\beta_1, \beta_2 \rangle=4np$ and $20 D_3-y^2=\det\langle\beta_1, \beta_3 \rangle=4mp$ for some integers $n\ge 2$ and $m\ge 2$ by Proposition \ref{prop:sublattice4np} and Proposition~\ref{prop:sublat4p}. Thus, $x= 2r$ and $y= 2s$ for some integers $r, s$ with $0\le r, s \le 5$. Since $D_2 \le D_3 \le \frac{3p}{5} + 5$ by Theorem~\ref{thm:D3notoverFp} and $p \geq 113$, we must have $m, n \in \{2, 3\}$.  

We show that in fact $n=2$ and $m=3$: Let $\{\beta_1,\beta^*_2, \beta^*_3\}$ be the Gram-Schmidt orthogonalization of the basis $\{\beta_1, \beta_2, \beta_3\}$. If $n=3$, 
we have that 
\begin{equation*}
D_2 \ge  \norm{\beta_2^*}^2= \frac{\det\langle\beta_1, \beta_2\rangle}{\norm{\beta_1}^2}= \frac{12p}{20}= \frac{3p}{5},
\end{equation*}
as well as
 \begin{equation*}
 \norm{\beta_3^*}^2= \frac{\det\curlyO^T}{\det\langle  \beta_1, \beta_2\rangle}= \frac{4p^2}{12p}=\frac{p}{3}
 \end{equation*}  
 and hence by equation \eqref{eq:d3bound1} and since $p \ge 113$, we have the following contradiction:
\[
D_3 \leq \norm{\beta^*_3}^2  + \frac{1}{4}\left(\norm{\beta_1}^2 + \norm{\beta_2^*}^2\right) = \frac{p}{3}+   \frac{1}{4}\left(\frac{3p}{5} + 20\right) < \frac{3p}{5} \leq D_2.
\]
Thus, $n=2$ and $D_2= \frac{2p+r^2}{5}$. Moreover, since $p\ge 113$, and $D_3 = \frac{mp+s^2}{5}$ with $m \in \{2,3\}$ and $0 \leq s \leq 5$,
\[D_3 \ge \norm{\beta_3^*}^2= \frac{\det\curlyO^T }{ \det\langle  \beta_1, \beta_2\rangle}= \frac{4p^2}{8p}=\frac{p}{2} > \frac{2p+25}{5},\]
and we must have $m=3$.

To complete the proof, it remains to narrow down further the values that the integers $z$, $r$ and $s$ can take. Using equation
\begin{equation*}
4p^2= \det\curlyO^T= \frac{1}{5}(24 p^2 -4 r^2 s^2 + 40 rsz  -100z^2),
\end{equation*}
we obtain that $z= \frac{rs \pm p}{5}$.
A case-by-case analysis considering all possible integer values for $r$ and $s$ in the range $[0,5]$ and checking the condition that 
\begin{equation*}
4p|\det\langle \beta_2, \beta_3\rangle= D_1D_3-z^2
\end{equation*}
for  $p\equiv13, 17 \pmod{20}$ yields the following  possibilities: $(r,s) \in \{(3, 1), (3,4)\}$ if $p\equiv 13 \pmod{20}$, and $(r,s) \in \{(1, 3), (1,2)\}$ if $p\equiv 17 \pmod{20}$. 

Finally, we show that $(r,s)\ne (3,4)$ when $p\equiv 13 \pmod{20}$. Indeed, if $(r,s) = (3,4)$, we have $D_2= \frac{2p+9}{5}$ and $z= \frac{12\pm p}{5}$. Since $5\nmid (12- p)$, it follows that $z= \frac{p+12}{5}$. However, in this case we have $|z| > \frac{D_2}{2}$, contradicting \Cref{lem:size_reduced}. Thus, $(r,s)=(3, 1)$, and since $5\nmid (3+p)$, we have $z= \frac{3-p}{5}$.

We can apply a similar argument for the case $p\equiv 17 \pmod{20}$ to obtain that  $(r,s)= (1,2 )$ and $z= \frac{2-p}{5}$.
\end{proof}

\subsection{The third successive minimum of Gross lattices when $j(E)\in \Fp$}\label {subsec:D3inFp}

In this remaining subsection, we prove the second part of our main result that characterizes the field of definition of the $j$-invariant of a supersingular elliptic curve by the value of the third successive minimum of its Gross lattice. When $p \equiv 3 \pmod{4}$, we also give necessary and sufficient conditions on $D_3$ for a maximal order in $B_p$ to be maximally embedded by $\mathbb{Z}[\frac{1+\sqrt{-p}}{2}]$, which is the second main result of this subsection.

\begin{theorem}\label{thm:D3overFppart1}
Let $p$ be a prime and $E$ be a supersingular elliptic curve defined over $\overline{\mathbb{F}}_p$ 
and $D_3$ be the third successive minimum of its Gross lattice. When $p \geq 7$, we have $j(E) \in \Fp\backslash\{0\}$ if and only if 
\begin{equation*}
p \le D_3  \le \frac{8p}{7} + \frac{7}{4}.
\end{equation*}
Furthermore, if $p \neq 3$ and $j(E) = 0$, then $D_3 = \frac{4p+1}{3},$ and finally if $p = 3$, we have $D_3 = 4$.    
\end{theorem}

\begin{remark}
We note that in fact this theorem covers all cases, as the unique supersingular elliptic curve defined over $\overline{\mathbb{F}}_p$ when $p = 2, 3, 5$ has $j$-invariant equal to $0$. However for these small primes we have $D_3 = \frac{4p+1}{3} < \frac{8}{7}p + \frac{7}{4}$. When $p = 7$, the only remaining case where $\frac{4p+1}{3} < \frac{8p}{7} + \frac{7}{4}$, the unique supersingular elliptic curve defined over $\overline{\mathbb{F}}_7$ has $j$-invariant equal to $6$ and $D_3 = 7 \leq \frac{8p}{7} + \frac{7}{4}$.
\end{remark}

\begin{proof}
The case of $j(E) = 0$ is handled in Proposition \ref{prop:0} for $p \neq 3$ and Remark \ref{rem:1728} for $p = 3$. 
In addition, the lower bound on $D_3$ follows from \Cref{lem:EnotinFp_D3lessthanp}.

To complete the proof, it remains to show that if $j(E)\in \Fp\backslash \{0\}$, then $D_3 \le \frac{8}{7}p + \frac{7}{4}$. As in the proof of Theorem \ref{thm:D3notoverFp}, we may size-reduce the successive minimal basis $\{\beta_1,\beta_2,\beta_3\}$ to obtain the basis $\{\beta_1,\beta_2,\beta'_3\}$. This size-reduced basis has the same Gram-Schmidt orthogonalization $\{\beta_1,\beta^*_2,\beta^*_3\}$ as our original successive minimal basis, and arguing in the same manner as we did to obtain equation \eqref{eq:d3bound}, we obtain the inequality
\begin{equation}\label{eq:d3boundFp}
D_3 \leq \norm{\beta^*_3}^2  + \frac{1}{4}\left(\norm{\beta_1}^2 + \norm{\beta^*_2}^2\right)= p+ \frac{1}{4}\left(X + \frac{4p}{X}\right),
\end{equation}
where as before $X=\max\{\norm{\beta^*_2}^2, \norm{\beta_1}^2 \}$, but this time $X \ge \sqrt{4p}$.
   
We consider the following three cases:
\begin{description}[leftmargin=*,itemindent=!]
    \item[Case 1] Let $X > \frac{p}{3}$ and $p>48$. Then by \Cref{lem:upperboundD1} with $n=1$,  we have that $\norm{\beta_1}^2< \frac{p}{3}$. Thus $\norm{\beta_2^*}^2= X > \frac{p}{3}$, and 
$$\norm{\beta_1}^2 = \frac{\det(\langle \beta_1, \beta_2 \rangle)}{\norm{\beta_2^*}^2}< \frac{4p}{p/3} = 12.$$
Since $\norm{\beta_1}^2 \equiv 0, 3 \pmod 4$, we must have $\norm{\beta_1}^2$ belonging to $\{3, 4, 7, 8, 11\}$. 

If $\norm{\beta_1}^2=3$, by Corollary \ref{cor:0and1728} $j(E)= 0$, which we have excluded from consideration. Similarly, if $\norm{\beta_1}^2=4$, again using Corollary \ref{cor:0and1728}, we conclude that $j(E)= 1728$, and by  Proposition~\ref{prop:1728}, this curve has $D_3= p+1 < \frac{8}{7}p + \frac{7}{4} $.

Finally if $\norm{\beta_1}^2 \in \{7, 8, 11\}$, then 
$$\norm{\beta_2^*}^2 = \frac{\det(\langle \beta_1, \beta_2 \rangle)}{\norm{\beta_1}^2}\le \frac{4p}{7}.$$
Applying Lemma~\ref{lem:max} with $a = 4p$ and $b = \frac{4p}{7}$ to equation \eqref{eq:d3boundFp}, and noting that $\sqrt{4p} < \frac{p}{3}<X$ since $p>48$, we have
\[
D_3 \le  p+ \frac{1}{4}\left(\frac{4p}{7} + \frac{4p}{4p/7}\right) = \frac{8}{7}p + \frac{7}{4}.
\]

\item[Case 2] Let $X \le \frac{p}{3}$ and $p > 48$. This time applying Lemma~\ref{lem:max} with $a = 4p$ and $b = \frac{p}{3}$ to equation \eqref{eq:d3boundFp}, and remembering that $X \ge \sqrt{4p}$, we get that 
\[
D_3 \le  p+ \frac{1}{4}\left(\frac{p}{3} + \frac{4p}{p/3}\right) \le \frac{8}{7}p + \frac{7}{4},
\]
since $p > 48$.

    \item[Case 3] We finally handle the case of $p \le 48$ with a finite computation. For each $p \le 48$, we compute a basis for each maximal order in $B_p$, and then the Gram matrix for a successive minimal basis of the associated Gross lattice. Rejecting all cases where $D_1 = 3$, which correspond to elliptic curves with $j$-invariant equal to $0$, and all cases where $D_3<p$, which correspond to elliptic curves with $j$-invariant in $\mathbb{F}_{p^2}\setminus \mathbb{F}_p$, we can verify that $D_3 \le \frac{8}{7}p + \frac{7}{4}$ in all remaining cases. Again, the computations are available on GitHub\footnote{\url{https://github.com/gkorpal/minimal-gross/blob/main/scripts/finite_cases/ref/cases_1_100.txt}}.
\end{description}
\end{proof}

As for our bound in the case of $j(E)$ in $\mathbb{F}_{p^2} \setminus \mathbb{F}_p$, the upper bound $D_3 \leq \frac{8}{7}p + \frac{7}{4}$ for $j(E) \in \mathbb{F}_p \setminus \{0\}$ given in Theorem~\ref{thm:D3overFppart1} is also quite tight. Indeed, we will see that when it is supersingular, the elliptic curve with $j$-invariant $15^3$ has Gross lattice with third successive minimum $\frac{8}{7}p + \frac{1}{7} \leq D_3$ by Proposition~\ref{prop:jminus15cube} below. 

We now turn our attention to demonstrating how the value of the third successive minimum of the Gross lattice of a supersingular elliptic curve with $j(E) \in \mathbb{F}_p$ for $p \equiv 3 \pmod{4}$ and $j(E) \neq 1728$ allows us to determine the arithmetic endomorphism ring of the elliptic curve. As a reminder, if $j(E) \neq 1728$, the ring of endomorphisms defined over $\mathbb{F}_p$ of a supersingular elliptic curve defined over $\mathbb{F}_p$ is determined by the $j$-invariant of $E$ (and not the $\mathbb{F}_p$-isomorphism class of the curve), and is isomorphic to $\mathbb{Z}[\sqrt{-p}]$ or $\mathbb{Z}[\frac{1+\sqrt{-p}}{2}]$. Furthermore, the geometric endomorphism ring of $E$, $\End(E)$, is maximally embedded by the arithmetic endomorphism ring. To obtain our main result we will need two propositions:

\begin{proposition}\label{prop:D3overFppart2}
Let $p \neq 3$ be a prime congruent to $3 \pmod{4}$, and $E$ be a supersingular elliptic curve defined over $\overline{\mathbb{F}}_p$ with $j(E)\in \Fp\backslash\{1728\}$. In this case, if $\End(E)$ is maximally embedded by $\mathbb{Z}[\frac{1+\sqrt{-p}}{2}]$, then $D_3= p$ and there exists a normalized successive minimal basis $\{\beta_1,\beta_2,\beta_3\}$ of $\curlyO^T$ with Gram matrix given by
\begin{equation*}
G_{\{\beta_1,\beta_2,\beta_3\}} = \begin{pmatrix}
D_1 & t_{12} & 0 \\
t_{12} & D_2 & 0 \\
0 & 0 & p
\end{pmatrix}
\end{equation*} with $t_{12} \in \mathbb{Z}$, $0\le t_{12} \leq \frac{D_1}{2}$ and $D_1 D_2 - t_{12}^2= 4p$. 
\end{proposition}

\begin{proof}
By \cite[Theorem 2]{Ibukiyama} and using his notation, if $\End(E)$ is maximally embedded by $\mathbb{Z}[\frac{1+\sqrt{-p}}{2}]$, we have $\End(E) \cong \curlyO'(q):= \mathbb{Z}\langle 1, \frac{1+\mj}{2}, \mi, \frac{r'\mi -\mk}{2q}\rangle$, where $q\equiv 3 \pmod 8$, $\left(\frac{p}{q} \right)=-1$, $\mi^2=-q, \mj^2=-p,$ and $(r')^2+p \equiv 0 \pmod {4q}$. Using this basis we can compute $\curlyO^T$ explicitly:
$$\curlyO'(q)^T= \langle 2\mi, \frac{r'\mi -\mk}{q}, \mj\rangle.$$
Because $(2\mi, \mj)= \left(\frac{r'\mi -\mk}{q}, \mj\right)=0$, $\mj$ is orthogonal to the sublattice $\Lambda'=\langle 2\mi, \frac{r'\mi -\mk}{q}\rangle$ and $\curlyO'(q)^T= \Lambda'\oplus \mathbb{Z} \mj$. Therefore we have $\det(\curlyO'(q)^T)= p \det(\Lambda')$ and hence $\det(\Lambda')= 4p$. Moreover, any successive minimal basis of $\curlyO'(q)^T$ must have one element contained in $\mathbb{Z} \mj$ and two others contained in $\Lambda'$.  
Therefore, $\norm{\mj}^2= p$ must be one of the successive minima of $\curlyO'(q)^T$. We show that $D_3=p$.

If $D_3\ne p$, then $D_1=p$ or $D_2=p$. $D_1 = p$ is impossible by the remark immediately following Theorem \ref{thm:kaneko} since we assume that $p \geq 7$. 

If $D_2= p$, then $ D_1 \le D_2=p \le D_3$ and the norms $D_1, D_3$ must be attained in $\Lambda'$. Thus,  $4p=\det(\Lambda') \le D_1 D_3 \le \frac{4}{3}\det(\Lambda')= \frac{16}{3}p$ using the Hermite bound (equation \eqref{eq:hermite}) for $\Lambda'$. This implies that $D_1\le \frac{16}{3}$ since $D_3\ge p$. 
Therefore we must have $D_1 \in \{3, 4\}$ and therefore $j(E) \in \{0, 1728\}$ by Corollary \ref{cor:0and1728}, or rather $j(E)=0$ since we exclude the case of $j(E) = 1728$ in this proposition. However, the endomorphism ring of the curve with $j$-invariant $0$ is isomorphic to the maximal order of the form $\curlyO(3)$ when $p \geq 5$ by \cite[Lemma 1.5]{Ibukiyama}, which is not maximally embedded by $\mathbb{Z}[\frac{1+\sqrt{-p}}{2}]$. Thus, we again have a contradiction, and $D_3=p$ is the third successive minimum of $\curlyO^T$. The form of the Gram matrix in the second statement then follows either by explicit computation or using Theorem \ref{thm:GLgrammatrix} for the entries that are equal to $0$ and Lemma \ref{lem:Gram_matrix} to obtain the bound on $t_{12}$, as well as Proposition \ref{prop:lalbet12}.
\end{proof}

We now turn our attention to the case of $j(E) = 1728$. In this case, if $p \equiv 3 \pmod{4}$ there are two $\mathbb{F}_p$-isomorphism classes of supersingular elliptic curves with $j(E) = 1728$ defined over $\mathbb{F}_p$, one of which has arithmetic endomorphism ring isomorphic to $\mathbb{Z}[\sqrt{-p}]$ and the other isomorphic to $\mathbb{Z}[\frac{1+\sqrt{-p}}{2}]$, and therefore the geometric endomorphism ring $\End(E)$ is maximally embedded by $\mathbb{Z}[\frac{1+\sqrt{-p}}{2}]$.

\begin{proposition}\label{prop:D_3_1728} 
Let $p\ge 7$, $E$ be a supersingular elliptic curve defined over $\overline{\mathbb{F}}_p$ and let $D_1 \leq D_2 \leq D_3$ be the successive minima of its Gross lattice. Then
    $$j(E)=1728\iff D_2=p \iff D_3=p+1.$$
\end{proposition}

\begin{proof}
We use the same notation as in Proposition~\ref{prop:1728}, and note that by Proposition~\ref{prop:1728}, if $j(E)= 1728$, then $D_2= p$ and $D_3=p+1$. Throughout, as usual $\{\beta_1,\beta_2,\beta_3\}$ is a successive minimal basis for the Gross lattice of $E$.

Now, assume that $D_2=p$. By Proposition \ref{prop:lalbet12}, we have that 
\begin{equation}\label{eq:1728}
\det(\langle \beta_1, \beta_2\rangle) = 4p= D_1 p - t_{12}^2,
\end{equation}
and by Theorem \ref{thm:GLgrammatrix}, $t_{12} = 0$, thus $D_1 = 4$ and $j(E)=1728$ by Corollary \ref{cor:0and1728}. By Proposition~\ref{prop:1728}, we also have that $D_3= p+1$.

Thus, to complete the proof of this corollary, it is sufficient to show that if $D_3= p+1$, then $j(E) = 1728$. If $D_3= p+1$, then by Theorem~\ref{thm:D3overFppart1}, $j(E) \in \Fp$. Now the lattices $\langle \beta_1, \beta_2\rangle$ and $\langle \beta_1, \beta_3\rangle$, respectively, have determinant $4p$ and $4np$ for some positive integer $n$, respectively, by Propositions \ref{prop:sublattice4np} and~\ref{prop:lalbet12}. Thus, 
\begin{equation}\label{eq:1st}
4p= D_1 D_2 - t_{12}^2,
\end{equation}
and
\begin{equation}\label{eq:1bst}
4np= D_1(p+1) - t_{13}^2,
\end{equation}
where again by Lemma \ref{lem:Gram_matrix} $t_{12},t_{13} \in \mathbb{Z}$ and $0\leq t_{12}, t_{13} \le \frac{D_1}{2}$.
Equation \eqref{eq:1bst} forces
\begin{equation}\label{eq:2nd}
p|(D_1 -t_{13}^2),
\end{equation}
and since $0\le t_{13} \le \frac{D_1}{2} \le \frac{2\sqrt{p}}{\sqrt{3}}$ by Theorem \ref{thm:kaneko}, it follows that
\begin{equation}\label{eq:3rdb}
-\frac{4p}{3}\le -\frac{D_1^2}{4}<D_1-t_{13}^2\le D_1 \le \frac{4\sqrt{p}}{\sqrt{3}} .
\end{equation}
Using equations \eqref{eq:2nd} and \eqref{eq:3rdb} and since $p\ge 7$, we must have $D_1-t_{13}^2 = -p$ or $D_1-t_{13}^2=0$. If $D_1-t_{13}^2 = -p$, by equation \eqref{eq:1bst} we have $D_1= 4n +1$, which contradicts the fact that $D_1 \equiv 0, 3 \pmod 4$. Thus, we must have $D_1= t_{13}^2$.

Carrying on, using Lemma \ref{lem:Gram_matrix} we also have 
\begin{equation}\label{eq:4th}
4p^2= \det(\curlyO^T)= (p+1)(D_1 D_2- t_{12}^2) - D_2y^2- D_1 t_{23}^2 + 2t_{12}t_{13}t_{23}. 
\end{equation}
Since $D_1 D_2-t_{12}^2 =4p$ and $D_1 = t_{13}^2$, equation \eqref{eq:4th} can be simplified:
\[4p^2= 4p^2 - (t_{12}-t_{13}t_{23})^2, \]
which implies $t_{12}=t_{13}t_{23}$. Hence, by equation \eqref{eq:1st}, we have 
\begin{equation*}
4p= t_{13}^2 D_2 - t_{13}^2 t_{23}^2= t_{13}^2( D_2-t_{23}^2).
\end{equation*}
and $t_{13}^2|4$ since $t_{13}^2= D_1 \le  \frac{4\sqrt{p}}{\sqrt{3}} <p$ for  $p\ge 7$. Since $t_{13}^2 =D_1\ge 3$, we must  have $D_1= t_{13}^2=4$, and $j(E) = 1728$ as claimed. 
 \end{proof}

Putting these two propositions together we obtain:

\begin{theorem}\label{thm:OqvsOqprime}
Let $p \neq 3$ be a prime, $E$ be a supersingular elliptic curve with $j$-invariant in $\Fp$, and $D_3$ be the third successive minimum of its Gross lattice. Then the following three statements are equivalent:
\begin{enumerate}
\item \label{state:1} $\End(E)$ is maximally embedded by $\mathbb{Z}[\frac{1+\sqrt{-p}}{2}]$ (and hence $p \equiv 3 \pmod{4}$);
\item \label{state:2} $D_3= p$ or $j(E)=1728$;  
\item \label{state:3} $D_3 \in \{p, p+1\} $. 
\end{enumerate}
\end{theorem}

\begin{proof}
Statements (\ref{state:2}) and (\ref{state:3}) are equivalent by \Cref{prop:D_3_1728}.
    
By Propositions~\ref{prop:0}, \ref{prop:D3overFppart2} and \ref{prop:D_3_1728}, to prove that (\ref{state:1}) and (\ref{state:2}) are equivalent, it suffices to show that if $j(E) \notin \{0, 1728\}$ and $D_3= p$, then $\End(E)$ is maximally embedded by $\mathbb{Z}[\frac{1+\sqrt{-p}}{2}]$.

As usual, let $\curlyO$ be a maximal order in $B_p$ such that $\curlyO\cong \End(E)$, and suppose that $D_3=p$ for the Gross lattice $\curlyO^T$. Then there is $\beta \in \curlyO^T$ with $\norm{\beta}^2 = p$ attaining this successive minimum, and by Proposition \ref{prop:endo_degd}, this element corresponds to an embedding of the imaginary quadratic order of discriminant $-p$, which is none other than $\mathbb{Z}[\frac{1+\sqrt{-p}}{2}]$, into $\curlyO$. Finally, this embedding must be optimal since $\beta$ must be primitive in $\curlyO^T$, since if $\beta = n \gamma$ for $\gamma \in \curlyO^T$ with $n \geq 2$, then $\norm{\gamma}^2 = \frac{p}{n^2}$ which is not an integer. This completes the proof.
\end{proof}

\section{Geometry of Gross lattices}\label{sec:geometry}
In this section, we present more details on the geometry of Gross lattices, by which we mean the norms of an ordered set of basis vectors (in our case, a successive minimal basis) and their pairwise inner product, which can be thought of as the ``angle" between two basis vectors. This information is exactly that contained in the Gram matrix of a normalized successive minimal basis of the Gross lattice.

We begin by investigating the existence of an orthogonal basis for the Gross lattice of a supersingular elliptic curve, and its well-roundedness in \Cref{subsec:mainorthorounded}. Following this, in \Cref{sec:uniqueness} we turn our attention to showing the uniqueness of the Gram matrix of a normalized successive minimal basis, considering the last few cases not covered by \cite[Corollary 3.15 and Section 4.1]{Goren-Love}. Then in \Cref{sec:algorithm} we investigate the structure of the Gram matrix of a successive minimal basis of the Gross lattice when $j(E) \in \mathbb{F}_p$, finding it is determined by the first two successive minima $D_1$ and $D_2$. This allows us to give an algorithm to compute the possible normalized Gram matrices for a supersingular elliptic curve with $j$-invariant in $\mathbb{F}_p$ given $p$ and the first successive minimum $D_1$ of its Gross lattice. Finally, in Section \ref{sec:special_cuves}, we use this algorithm to compute the normalized Gram matrix of the reduction modulo $p$ of each of the $13$ $\overline{\mathbb{Q}}$-isomorphism classes of elliptic curves that have complex multiplication by an imaginary quadratic order of class number $1$ when this reduction is supersingular and this imaginary quadratic order is the order of greatest discriminant that injects into the endomorphism ring of the reduction.

\subsection{Orthogonality and well-roundedness of Gross lattices}\label{subsec:mainorthorounded}
A lattice is \textbf{orthogonal} if it has an orthogonal basis, i.e., every pair of distinct vectors in this basis has inner product equal to zero. We note that in this case, the orthogonal basis, when its elements are ordered by increasing norm, is also a successive minimal basis. A lattice of rank $n$ is \textbf{well-rounded} if it has $n$ linearly-independent shortest vectors, i.e., all of its $n$ successive minima are equal. A natural question when considering the geometry of a lattice is whether it is orthogonal or well-rounded, and we turn our attention to these questions for the Gross lattice of a supersingular elliptic curve now. 

\begin{theorem}\label{thm:orthogonal_WR}
Let $p$ be a prime and $E$ be a supersingular elliptic curve defined over $\overline{\mathbb{F}}_p$.

If $p=2$, the unique supersingular elliptic curve defined over $\overline{\mathbb{F}}_2$ has a well-rounded Gross lattice, but this lattice is not orthogonal.   
  
If $p\ge 3$, then the Gross lattice of a supersingular curve defined over $\overline{\mathbb{F}}_p$ is neither orthogonal nor well-rounded.
\end{theorem}

To prove Theorem \ref{thm:orthogonal_WR}, we give more precise results: The first about the equality of successive minima of the Gross lattice of a supersingular elliptic curve, and the second giving necessary and sufficient conditions for $t_{ij} = (\beta_i,\beta_j)$ to be equal to zero.

First, recall that as stated here in Theorem \ref{thm:GLinequality}, which is adapted from \cite[Lemmata 4.4 and 4.5]{Goren-Love}, we have that if $p \geq 11$, then $D_1 \neq D_2$ and if in addition $D_1 \geq 15$, then $D_2 \neq D_3$. As a consequence of our work we can remove the hypotheses on this result and determine exactly when we have equality of successive minima for the Gross lattice of a supersingular elliptic curve:

\begin{proposition}\label{prop:inequalityDi}
Let $p$ be a prime and $E$ be a supersingular elliptic curve defined over $\overline{\mathbb{F}}_p$ and $D_1 \leq D_2 \leq D_3$ be the successive minima of its Gross lattice. If $p \neq 2$, then $D_1\ne D_2$. Furthermore, if $j(E) \neq 0$, then $D_2 \neq D_3$.
\end{proposition}

\begin{proof}
    The cases where $2 \leq p \leq 11$ can be checked with a finite computation\footnote{\url{https://github.com/gkorpal/minimal-gross/blob/main/scripts/finite_cases/ref/cases_1_100.txt}}; we use that $D_1 = 3$ if and only if $j(E) = 0$ to check that $D_2=D_3$ only for those cases.

    Now if $p \geq 13$ and $D_1 < 15$, the imaginary quadratic order of discriminant $-D_1$ has class number $1$ by the results compiled in Table \ref{tab:CMorders}. Hence by Corollary \ref{cor:ssreduction}, $j(E) \in \mathbb{F}_p$, and by Proposition \ref{prop:sublattice4np} we have  $\det \langle \beta_2, \beta_3\rangle = 4 np$ for some positive integer $n$. First, we use that $4np = D_2^2- t_{23}^2 = (D_2- t_{23})(D_2 +t_{23})$ for $t_{23}= (\beta_2, \beta_3)$ and $|t_{23}| \le \frac{D_2}{2}$ by \Cref{lem:size_reduced}, thus it follows that $p|(D_2-t_{23})$ or $p|(D_2+t_{23})$. If $D_2 = D_3$, by \Cref{thm:D3overFppart1}, $p \le D_2=D_3 \le \frac{8}{7}p + \frac{7}{4}$ for $j(E)\ne 0$, and so for $p\ge 11$ we have
$$\frac{p}{2}\le \frac{D_2}{2}\le D_2  \pm t_{23} \le \frac{3D_2}{2} \le \frac{12}{7}p + \frac{21}{8} < 2p.$$ 
Therefore, we either have that $D_2-t_{23}=p$ and $D_2+t_{23}= 4n$, or that $D_2+t_{23}=p$ and $D_2-t_{23}= 4n$. In both cases, we have that $D_2= \frac{p+4n}{2}$, which is not an integer, leading to a contradiction. Finally if $j(E) = 0$ then $D_2 = D_3$ by Proposition \ref{prop:0}.
\end{proof}

We can also determine exactly when the elements of a successive minimal basis are orthogonal:

\begin{proposition}\label{prop:tijzero}
    Let $p$ be a prime and $E$ be a supersingular elliptic curve defined over $\overline{\mathbb{F}}_p$, $\mathcal{O}$ be a maximal order in $B_p$ isomorphic to its endomorphism ring, $\{\beta_1,\beta_2,\beta_3\}$ a successive minimal basis of $\mathcal{O}^T$, and $t_{ij} = \frac{1}{2}\trd(\beta_i\overline{\beta}_j)$. Then we have:
    \begin{enumerate}
    \item $t_{12} = t_{13}= 0$ if and only if $p = 3$, and in this case $t_{23} \neq 0$;
    \item $t_{12} = t_{23}=0$ if and only if $p \neq 3$ and $j(E) = 1728$, and in this case $t_{13} \neq 0$;
    \item $t_{13}= t_{23} = 0$ if and only if $j(E) \neq 1728$ and $\mathcal{O}$ is maximally embedded by $\mathbb{Z}[\frac{1+\sqrt{-p}}{2}]$, and in this case $t_{12} \neq 0$;
    \item in all other cases, $t_{ij} \neq 0$ for any $i$ and $j$.   
    \end{enumerate}
\end{proposition}

\begin{proof}
    By \Cref{rem:uniquej}, we have that $D_1 < p$ except if $p = 2, 3$ or $5$. For each of these 3 primes, the unique supersingular curve over $\overline{\mathbb{F}}_p$ has $j(E)=0$ and hence $D_1 = 3$.

    When $p=2$, an argument entirely analogous to that given in Remark \ref{rem:1728} for $p = 3$ shows that the Gram matrix of a normalized successive minimal basis for the Gross lattice must be
    \begin{equation*}
    \begin{pmatrix}
        3 & 1 & 1 \\
        1 & 3 & -1 \\ 
        1 & -1 & 3
    \end{pmatrix},
    \end{equation*}
    and therefore $t_{ij} \neq 0$ for any $i,j$. If $p = 3$, using Remark \ref{rem:1728} we obtain that $t_{12} = t_{13} = 0$ but $t_{23} \neq 0$. 
    
    If $p \geq 5$, then $D_1 < p$ and in particular $D_1 \neq p$. For each $i,j$, we have that $D_iD_j - t_{ij}^2$ is a multiple of $p$ by Proposition \ref{prop:sublattice4np}, and therefore if $t_{ij} = 0$ then $D_i$ or $D_j$ is a multiple of $p$. By Theorems \ref{thm:D3notoverFp} and \ref{thm:D3overFppart1}, $D_3< 2p$, and so in fact if $t_{ij} = 0$ then $D_i$ or $D_j$ is equal to $p$, and here the ``or" is exclusive as by Proposition \ref{prop:inequalityDi}, we never have $D_2 = D_3 = p.$ We conclude thus that it is impossible for $t_{12} = t_{13} = t_{23} =0$, or for $t_{12} = t_{13} = 0$ when $p \geq 5$, and that if any $t_{ij} =0$ then $j(E) \in \mathbb{F}_p$, since if $j(E) \in \mathbb{F}_{p^2} \setminus \mathbb{F}_p$ then $D_3 < p$.
    
    By Proposition \ref{prop:D_3_1728}, $D_2 = p$ if and only if $j(E) = 1728$. In that case, by Theorem \ref{thm:GLgrammatrix}, we have $t_{12} =t_{23}=0$. 
    
    Finally, $D_3 = p$ if and only if $\mathcal{O}$ is maximally embedded by $\mathbb{Z}[\frac{1+\sqrt{-p}}{2}]$ but $j(E) \neq 1728$. In that case, again by Theorem \ref{thm:GLgrammatrix}, $t_{13} = t_{23}= 0$.
\end{proof}

\begin{remark}\label{rem:choice}
    By Theorem \ref{thm:GLgrammatrix}, the Gram matrix of a successive minimal basis of $\mathcal{O}^T$ is equal to the Gram matrix of any other successive minimal basis, up to the sign changes on the values $t_{ij}$ induced by sending $\beta_i$ to $-\beta_i$ for $1 \leq i \leq 3.$ Hence if all $t_{ij} \neq 0$, fixing the sign of any two of them fixes the sign of the third, and in turn the value of the normalized Gram matrix. However as we saw there are cases where two of the values $t_{ij}$ can be zero and since $t_{12} = t_{13} = 0$ if and only if $p = 3$, fixing $t_{12}$ and $t_{13}$ to be nonnegative fixes the Gram matrix in all cases except one, as we show in Theorem \ref{thm:Gram_unique} below. This explains our choice of convention that $(\beta_1,\beta_2)$ and $(\beta_1,\beta_3)$ be nonnegative in a normalized successive minimal basis, rather than any other choice.
\end{remark}

\subsection{Uniqueness of the normalized Gross lattices }\label{sec:uniqueness}

By Theorem \ref{thm:GLgrammatrix}, the Gram matrix of a normalized successive minimal basis is unique, i.e, independent of the choice of basis, if $p$ is odd and $D_1 \geq 8$. We now remove these hypotheses to show that this Gram matrix is unique except in the case of $p = 3$:

\begin{theorem}\label{thm:Gram_unique}
Let $p \neq 3$ and $E$ be a supersingular elliptic curve defined over $\overline{\mathbb{F}}_p$. 
Then all normalized successive minimal bases of the Gross lattice of $E$ have the same Gram matrix $G_{\curlyO^T}$, whose form is given in \Cref{lem:Gram_matrix}. 
\end{theorem}

\begin{proof}
The case of $p = 2$ is handled in the proof of Proposition \ref{prop:tijzero}, hence we may assume from now that $p$ is odd, and therefore by Theorem \ref{thm:GLgrammatrix}, given elements $\beta_i$ and $\beta_j$ of $\mathcal{O}^T$ attaining the successive minima $D_i$ and $D_j$ for $i \neq j$, the absolute value of the inner product $\left|(\beta_i,\beta_j)\right|$ is determined by the values $D_i$ and $D_j$ as long as $\min \{ D_i,D_j\} \leq p$. In that case, if the product $t_{12}t_{13}t_{23} \neq 0$, fixing the sign of $t_{12}$ and $t_{13}$ fixes the sign of $t_{23}$, and hence we obtain a unique normalized Gram matrix. As shown in Proposition \ref{prop:tijzero}, if $p \neq 3$ and $t_{12}t_{13}t_{23} = 0$ then $t_{23} =0$, and therefore again fixing the sign of $t_{12}$ and $t_{13}$ yields a unique normalized Gram matrix.

Hence, to complete the proof, it suffices to handle the cases of $D_1 = 3,4$ and $7$. We remark that the hypothesis that $D_1 \geq 8$ is only used in \cite[Section 4.1]{Goren-Love} to ensure that $D_2 \leq p$, which in turns allows the use of Corollary 3.15 of \emph{loc.\ cit.}\ to ensure that the absolute value of $t_{23}$ is determined by the values $D_2$ and $D_3$. Therefore to handle the case of $D_1=7$, it suffices to show that in this case as well $D_2 \leq p$, and the argument of Section 4.1 of \emph{loc.\ cit.}\ applies. But this is true, as if $D_1 = 7$, then $j(E) \in \mathbb{F}_p$, and using the Hermite bound \eqref{eq:hermite} as well as the fact that $\det \langle \beta_1,\beta_2\rangle = 4p$ by Proposition \ref{prop:lalbet12}, we have $D_2 \leq \frac{16p}{21}$.

Finally we may consider the cases $D_1=3$, when $j(E) =0$, and $D_1 = 4,$ when $j(E) = 1728$. Since we have handled the case of $p =2$ and we exclude the case of $p =3$, we may assume that $p \geq 5$ which ensures that $D_1 \leq p$ and we can apply Theorem \ref{thm:GLgrammatrix} to determine the absolute value of $t_{12}$ and $t_{13}$.

Now consider the case $j(E)=0$ and $p \neq 3$. By \Cref{prop:0}, the Gross lattice $\curlyO^T$ of $E$ has $D_1= 3$ and $D_2= D_3= \frac{4p+1}{3},$ and by \ref{thm:GLgrammatrix} we have $t_{12}= t_{13} =1$ for any normalized successive minimal basis. Computing the determinant of $\curlyO^T$ from its entries, we obtain
\begin{equation*}
\det(\curlyO^T)= 4p^2=\frac{1}{3} (16p^2 - (3t_{23}-1)^2),
\end{equation*}
which implies that $t_{23}= \frac{\pm 2p +1}{3}$. Because we have that $|t_{23}| \leq \frac{D_2}{2}= \frac{4p+1}{6}$, we must have $t_{23}= \frac{- 2p +1}{3}$, and this value is unique. 

The case $j(E)= 1728$ and $p \neq 3$ can be obtained using a similar argument as in the case of $j(E)=0$ using \Cref{prop:1728}. 
\end{proof}

\subsection{The possible normalized Gram matrices of supersingular curves with $j(E) \in \mathbb{F}_p$}\label{sec:algorithm}

In this section we begin by giving a result showing that if $j(E) \in \mathbb{F}_p$, the normalized Gram matrix of $\mathcal{O}^T$ is determined completely by the values of $D_1$ and $D_2$. A much simplified proof was shared with the authors by Jonathan Love, and we present it in Appendix \ref{appendix}. Using this result and those of \Cref{sec:mainsection}, we can then give an algorithm (\Cref{alg:gramgrossD1}) to compute the possible values of a normalized Gram matrix for a supersingular elliptic curve with $j$-invariant in $\mathbb{F}_p$ given $p$ and the value $D_1$ of its first successive minimum. 

\begin{theorem}\label{thm:4types}
Let $p$ be a prime and $E$ be a supersingular elliptic curve defined over $\overline{\mathbb{F}}_p$ with $j(E)\in \Fp$, and $\curlyO$ be a maximal order of $B_p$ isomorphic to $\End(E)$. Then the normalized Gram matrix $G_{\mathcal{O}^T}$ of the Gross lattice of $\curlyO$ has one of the following forms:

\begin{itemize}
   \item \textbf{Type 1,} $D_1 \equiv D_2 \equiv 0 \pmod 4$:
\[G_{\mathcal{O}^T}=
\begin{pmatrix}
D_1 & t_{12} & 0 \\
t_{12} & D_2 & 0 \\
0 & 0 & p
\end{pmatrix}.
\]

\item \textbf{Type 2,} $D_1 \equiv 0 \pmod{4}$, and $D_2 \equiv 3 \pmod 4$:
\[G_{\mathcal{O}^T}=
\begin{pmatrix} 
D_1 & t_{12} & \tfrac{D_1}{2} \\
t_{12} & D_2 & \tfrac{t_{12}}{2} \\
\tfrac{D_1}{2} & \tfrac{t_{12}}{2} & p + \tfrac{D_1}{4}
\end{pmatrix}.
\]

\item \textbf{Type 3,} $D_1 \equiv 3 \pmod{4}$, and $D_2 \equiv 0 \pmod 4$:
\[G_{\mathcal{O}^T}=
\begin{pmatrix}
D_1 & t_{12} & \tfrac{t_{12}}{2} \\
t_{12} & D_2 & \tfrac{D_2}{2} \\
\tfrac{t_{12}}{2} & \tfrac{D_2}{2} & p + \tfrac{D_2}{4}
\end{pmatrix}.
\]

\item \textbf{Type 4,} $D_1 \equiv D_2 \equiv 3 \pmod 4$:
\[G_{\mathcal{O}^T}=
\begin{pmatrix}
D_1 & t_{12} & \tfrac{D_1 - t_{12}}{2} \\
t_{12} & D_2 & \tfrac{t_{12} - D_2}{2} \\
\tfrac{D_1 - t_{12}}{2} & \tfrac{t_{12} - D_2}{2} & p + \tfrac{D_1 - 2t_{12} + D_2}{4}
\end{pmatrix},
\]
\end{itemize}
where here as in the rest of the article, $D_1\le D_2$ are the first two successive minima of the Gross lattice of $E$, $D_i \equiv 0, 3 \pmod 4$ for $i = 1,2$, $t_{12} := \frac{1}{2}\trd(\beta_1\overline{\beta_2})$, $0\le t_{12} \le \frac{D_1}{2}$, and since $j(E) \in \mathbb{F}_p$, $D_1 D_2- t_{12}^2= 4p$.
\end{theorem}

\begin{remark}
Using the notation of \Cref{thm:ibukiyama}, we note that if $j(E) \neq 1728$, if $\curlyO=\curlyO'(q',r')$ then $G_{\mathcal{O}^T}$ is of Type 1, which agrees with Proposition \ref{prop:tijzero} stating that $t_{13} = t_{23} = 0$ if and only if $\mathcal{O}$ is embedded by $\mathbb{Z}[\frac{1+ \sqrt{-p}}{2}]$ and $j(E) \neq 1728$. Still in the case of $j(E) \neq 1728$, if $\curlyO=\curlyO(q,r)$ then $G_{\mathcal{O}^T}$ is of Type 2, 3 or 4. Finally, applying Proposition \ref{prop:tijzero}, if $j(E) = 1728$, then $G_{\mathcal{O}^T}$ is of Type 3 if $p = 3$ and of Type 2 otherwise.
\end{remark}

The result of \Cref{thm:4types} shows that the normalized Gram matrix of the Gross lattice of a supersingular elliptic curve with $j$-invariant in $\Fp$ is determined by $D_1$ and $D_2$ only, since these values determine $t_{12}$ uniquely if $p \neq 3$ by Theorem \ref{thm:Gram_unique} as well as when $p = 3$ by Remark \ref{rem:1728}. 
Applying this result, we can give an algorithm to compute the normalized Gram matrix of the Gross lattice of supersingular elliptic curves with $j$-invariant in $\Fp$, which is presented as Algorithm \ref{alg:gramgrossD1} below.

\begin{algorithm}[!ht] 
\caption{$\mathsf{GramGross}(p,D_1)$}\label{alg:gramgrossD1} 
    \begin{algorithmic}[1]
    \REQUIRE A prime $p$ and a positive integer $D_1 \equiv 0, 3 \pmod 4$ such that $D_1\leq \frac{4\sqrt{p}}{\sqrt{3}}$.
    
    \ENSURE All possible normalized Gram matrices for a supersingular elliptic curve with $j$-invariant in $\Fp$, and such that the first successive minimum is $D_1$. 
    
    \STATE $L \leftarrow \texttt{[]}$ \COMMENT{Initialize an empty list of Gram matrices}
    \STATE $L_{12} \leftarrow \texttt{[]}$ \COMMENT{Initialize an empty list of eligible $[t_{12},D_2]$ pairs.} 
    \STATE \label{alg:checksquarebegin}  $B = \lfloor D_1/2\rfloor$ 
    \FOR{$t_{12} = 0$ \TO $B$ \label{alg:for1}}
    \STATE $D_2 \leftarrow (4p+t_{12}^2)/D_1$ 
    \IF{$D_2\in \mathbb{Z}$ \AND $D_1\leq D_2$ \AND $D_2\equiv 0,3\pmod 4$}
    \STATE $L_{12}.\text{append}([t_{12},D_2])$
    \ENDIF
    \ENDFOR \label{alg:checksquareend}
    \IF{$L_{12} = \emptyset$}
    \RETURN $L$ \COMMENT{No Gram matrix.}
    \ENDIF
    
    \FOR{$[t_{12},D_2]\in L_{12}$ \label{alg:matrixbegin}}     
    \IF{$D_1 \equiv D_2 \equiv 0\pmod 4$}
    \STATE $t_{13} \leftarrow 0,  t_{23} \leftarrow 0, D_3 \leftarrow p$ 
    
    \ELSIF{$D_1 \equiv  0\pmod 4$ and $D_2 \equiv 3\pmod 4$}
    \STATE $t_{13} \leftarrow \frac{D_1}{2}, t_{23} \leftarrow \frac{t_{12}}{2}, D_3 \leftarrow p + \frac{D_1}{4}$
    
    \ELSIF{$D_1 \equiv  3\pmod 4$ and $D_2 \equiv 0\pmod 4$}
    \STATE $t_{13} \leftarrow \frac{t_{12}}{2}, t_{23} \leftarrow \frac{D_2}{2}, D_3 \leftarrow p + \frac{D_2}{4}$

    \ELSIF{$D_1 \equiv  3\pmod 4$ and $D_2 \equiv 3\pmod 4$}
    \STATE \label{alg:matrixend} $t_{13} \leftarrow \frac{D_1-t_{12}}{2}, t_{23} \leftarrow \frac{t_{12}-D_2}{2}, D_3 \leftarrow p + \frac{D_1-2t_{12}+D_2}{4}$
    
    \STATE $L.\text{append}(G)$
    \ENDIF
    \ENDFOR
    
    \RETURN $L$ 
    \end{algorithmic}
\end{algorithm}

We make two remarks about this algorithm: 
First, in the case where $p = 3$ and the unique isomorphism class of supersingular elliptic curves over $\overline{\mathbb{F}}_3$ does not admit a unique normalized Gram matrix, the algorithm returns the Gram matrix of a successive minimal basis such that $(\beta_2,\beta_3)>0$. 
Secondly, if $-4p$ is not a square modulo $D_1$, then no Gram matrix will be found. If the squares modulo $D_1$ are known (for example if $D_1$ is a product of distinct primes and quadratic reciprocity can be used, or if $D_1$ is small), performing this additional check before the beginning of the algorithm can avoid running it unnecessarily. 

\begin{proposition}\label{Correcness_runtime}
\Cref{alg:gramgrossD1} is correct and its run time is $O(D_1)$.
    
\end{proposition}

\begin{proof}
Lines \ref{alg:checksquarebegin} to \ref{alg:checksquareend} in the algorithm follow from \Cref{prop:lalbet12} and Lemma \ref{lem:size_reduced}, and Lines \ref{alg:matrixbegin} to \ref{alg:matrixend} from \Cref{thm:4types}. The algorithm is then correct.

As for the run time, the length of the two for loops (beginning on lines \ref{alg:for1} and \ref{alg:matrixbegin}) depend linearly on the size of $D_1$ and on the size of $L_{12}$, respectively, which are both $O(D_1)$.
\end{proof}

Finally, we note that we were not able to prove that each of the output matrices of \Cref{alg:gramgrossD1} corresponds to the Gram matrix of a normalized successive minimal basis for the Gross lattice of a supersingular elliptic curve with $j$-invariant in $\Fp$. However, we have found this to be true in our experiments for all primes $p$ up to $10^6+3$. 

\subsection{Gram matrices of Gross lattices of special curves}\label{sec:special_cuves}

We end this article by applying \Cref{alg:gramgrossD1} to compute the Gram matrix of a successive minimal basis for the supersingular reduction of each of the $13$ elliptic curves defined over $\overline{\mathbb{Q}}$ with complex multiplication by an imaginary quadratic order of class number $1$ listed in \Cref{tab:CMorders} when $p$ is large enough that this imaginary quadratic order is the order of largest discriminant that embeds in the endomorphism ring. (We may apply the algorithm as each of these reductions has $j(E) \in \mathbb{F}_p$ by Corollary \ref{cor:ssreduction}). To do this, in Table \ref{table:CMcurves2} we give for each curve $E$ with complex multiplication by an imaginary quadratic order of class number $1$ and discriminant $-d$, we give a sharp lower bound $N_E$ such that if $p \geq N_E$ and $E$ has supersingular reduction at $p$, then $D_1$, the first successive minimum of the reduction modulo $p$ of $E$, satisfies $D_1 = d$. We also compare this sharp, experimentally obtained bound $N_E$ to the bound $(d+1)^2/4$ obtained in \cite[Lemma 2.2.1]{Goren-Lauter}. In each case, our algorithm produced a single possible Gram matrix, which allows us to conclude that we have found the normalized Gram matrix of the unique supersingular elliptic curve defined over $\overline{\mathbb{F}}_p$ whose endomorphism ring is maximally embedded by this quadratic order. 

\begin{lemma}\label{lem:D1equalsdisc}
Let $E$ be the (unique up to $\overline{\mathbb{Q}}$-isomorphism) elliptic curve defined over $\overline{\mathbb{Q}}$ with complex multiplication by the imaginary quadratic order $O$ of discriminant $-d$, where $O$ is one of the orders listed in Table~\ref{tab:CMorders}. 
Then there exists a least positive prime $N_E$, which depends only on $j(E)$, such that for all primes $p$ of good supersingular reduction for $E$ such that $p\geq N_E$, the first successive minimum of the Gross lattice of $\widetilde{E}$ is equal to $d$, where $\widetilde{E}$ is the reduction modulo $p$ of $E$. 
\end{lemma}

\begin{proof}  
Since the endomorphism ring of $E$ injects into the endomorphism ring of $\widetilde{E},$ by Proposition \ref{prop:endo_degd} the Gross lattice of $\widetilde{E}$ contains an element of norm $d$. By \cite[Lemma 2.2.1]{Goren-Lauter}, which we note is proved using techniques different from those used by Gross-Zagier \cite{Gross-Zagier84}, if $p > (d+1)^2/4$, the maximal order $\curlyO \cong \End(\widetilde{E})$ does not embed any imaginary quadratic order of discriminant with absolute value smaller than $d$. As a consequence, using Proposition \ref{prop:endo_degd} again, the first successive minimum of the Gross lattice of $\widetilde{E}$ must be $d$ for $p$ large enough when $p$ is a prime of supersingular reduction for $E$.
\end{proof}

We note that Lemma \ref{lem:D1equalsdisc} defines $N_E$ to be the \emph{least} prime such that $p \geq N_E$ implies that the first successive minimum of Gross lattice of the supersingular reduction of $E$ modulo $p$ is equal to $d$. As a consequence, the value of the least prime strictly greater than $(d+1)^2/4$ is an upper bound for the value of $N_E$. For the $13$ elliptic curves with complex multiplication by an imaginary quadratic order of class number $1$, we compute the value of $N_E$ exactly in the following manner: For each elliptic curve $E$ and for each prime $p$ of supersingular reduction of the curve smaller than $(d+1)^2/4$, we apply \cite[Algorithms 1 and 2]{LOX20} to the curve $\widetilde{E}$ which is the reduction modulo $p$ of $E$ to establish directly the value $D_1$ of the first successive minimum of its Gross lattice\footnote{\url{https://github.com/gkorpal/minimal-gross/tree/main/scripts/NE_values}}. Our results are listed in \Cref{table:CMcurves2}, and they demonstrate that except if $j(E)\in \{0, 1728, 20^3\}$, we have that $N_E < (d+1)^2/4$.

\begin{table}[!ht]
\caption{$j$-invariant, first successive minimum $D_1$, and constant $N_E$.}
\label{table:CMcurves2}
  \begin{tabular}{|l||wc{2em}|c|c|c|c|c|c|c|c|}
  \hline 
  $j$-invariant & $0$ & $1728$ & $-15^3$ & $20^3$ & $-32^3$ & $2\cdot 30^3$ & $66^3$ &  $-96^3$ \\ \hline
  $D_1=d$ & $3$ & $4$ & $7$ & $8$ & $11$ & $12$ & $16$ & $19$ \\ \hline
  $(d+1)^2/4$ & $4$ & $6.25$ & $16$ & $20.25$ & $36$ & $42.25$ & $72.5$ & $100$ \\ \hline
  $N_E$ & $5$ & $7$ & $13$ &  $23$ & $29$ &  $41$ & $67$ & $79$ \\ \hline
  \end{tabular}
  
 \bigskip
 
  \begin{tabular}{|l||c|c|c|c|c|}
  \hline 
  $j$-invariant & $-3\cdot 160^3$ &  $255^3$ &   $-960^3$ & $-5280^3$ & $-640320^3$ \\ \hline
  $D_1=d$ & $27$ &  $28$ & $43$ & $67$ & $163$ \\ \hline
  $(d+1)^2/4$ & $196$ & $210.25$ & $484$ & $1156$ & $6724$ \\ \hline
  $N_E$ & $167$ &  $181$ & $433$ & $1103$  & $6481$ \\ \hline
  \end{tabular}
\end{table}

With the value of $N_E$ in hand for each of the elliptic curves under consideration, we can then use SageMath to compute the normalized Gram matrix of the reduction modulo $p$ of $E$ when $p$ is a prime of supersingular reduction, treating the prime $p$ as a variable. Here is a representative result we have obtained from applying Algorithm \ref{alg:gramgrossD1}\footnote{\url{https://github.com/gkorpal/minimal-gross/blob/main/scripts/CM_Gross/ref/cm_7.txt}} to the supersingular elliptic curves with $j$-invariant $-15^3$:

\begin{proposition}\label{prop:jminus15cube}
    If  $p\geq 13$, then the normalized Gram matrix of the Gross lattice of the supersingular elliptic curve $E$ defined over $\overline{\mathbb{F}}_p$ with $j(E)=-15^3$ is one of the following:
    \[G_{\curlyO^T} = \begin{pmatrix}
7 & 3 & 2 \\
3 & \frac{4p+9}{7} & -\frac{2p-6}{7}\\
2 & -\frac{2p-6}{7} & \frac{8p+4}{7}
\end{pmatrix}\quad \text{if } p\equiv 3\pmod {7},\]
 \[G_{\curlyO^T} = \begin{pmatrix}
7 & 1 & 3 \\
1 & \frac{4p+1}{7} & -\frac{2p-3}{7}\\
3 & -\frac{2p-3}{7} & \frac{8p+9}{7}
\end{pmatrix}\quad \text{if } p\equiv 5\pmod {7},\]
or
 \[G_{\curlyO^T} = \begin{pmatrix}
7 & 2 & 1 \\
2 & \frac{4(p+1)}{7} & \frac{2(p+1)}{7}\\
1 & \frac{2(p+1)}{7} & \frac{8p+1}{7}
\end{pmatrix}\quad \text{if } p\equiv 6\pmod {7}.\]
\end{proposition}

Finally, for each remaining supersingular $j$-invariant in $\mathbb{Q}$ corresponding to the elliptic curve defined over $\overline{\mathbb{Q}}$ with complex multiplication by a quadratic order $O$ of class number $1$, in Table \ref{table:CMcurves} we list first the congruence conditions on $p$ ensuring that $p$ is inert in $K = O \otimes_\mathbb{Z} \mathbb{Q}$ (and hence that the reduction modulo $p$ of the curve $E$ is supersingular), and then a link to a file containing the normalized Gram matrix of the Gross lattice of the supersingular reduction, depending on the congruence class of $p$.

\begin{table}[H]
\caption{Supersingular primes \cite[Theorem 10.7]{washington} and links to their Gram matrix in our repository on GitHub.}
\label{table:CMcurves}
\begin{tabular}{|wr{5em}|p{20em}|p{6em}|}
  \hline
  $j$-invariant & Supersingular primes $p\equiv a \pmod d$ & Link to data\\ \hline \hline
  $0$ &  $2\pmod 3$ & \href{https://github.com/gkorpal/minimal-gross/blob/main/scripts/CM_Gross/ref/cm_3.txt}{\texttt{cm\_3.txt}} \\ \hline
  $12^3$ & $3\pmod 4$ & \href{https://github.com/gkorpal/minimal-gross/blob/main/scripts/CM_Gross/ref/cm_4.txt}{\texttt{cm\_4.txt}} \\ \hline
  $20^3$  &  $5,7\pmod 8$ & \href{https://github.com/gkorpal/minimal-gross/blob/main/scripts/CM_Gross/ref/cm_8.txt}{\texttt{cm\_8.txt}}  \\ \hline
  %; \href{https://github.com/gkorpal/minimal-gross/blob/main/scripts/CM_Gross_finite_cases/ref/cases_8.txt}{\texttt{cases\_8.txt}}
  $2\cdot 30^3$ &  $5,11\pmod {12}$ & \href{https://github.com/gkorpal/minimal-gross/blob/main/scripts/CM_Gross/ref/cm_12.txt}{\texttt{cm\_12.txt}} \\ \hline
  $-15^3$  & $3,5,6\pmod 7$ & \href{https://github.com/gkorpal/minimal-gross/blob/main/scripts/CM_Gross/ref/cm_7.txt}{\texttt{cm\_7.txt}}  \\ \hline
  %; \href{https://github.com/gkorpal/minimal-gross/blob/main/scripts/CM_Gross_finite_cases/ref/cases_7.txt}{\texttt{cases\_7.txt}}
  $66^3$ & $3,7,11,15 \pmod {16}$ & \href{https://github.com/gkorpal/minimal-gross/blob/main/scripts/CM_Gross/ref/cm_16.txt}{\texttt{cm\_16.txt}}  \\ \hline
  $-32^3$ & $2,6,7,8,10 \pmod {11}$ & \href{https://github.com/gkorpal/minimal-gross/blob/main/scripts/CM_Gross/ref/cm_11.txt}{\texttt{cm\_11.txt}} \\ \hline
  %; \href{https://github.com/gkorpal/minimal-gross/blob/main/scripts/CM_Gross_finite_cases/ref/cases_11.txt}{\texttt{cases\_11.txt}}
  $255^3$  & $3,5,13,17,19,27\pmod {28}$ & \href{https://github.com/gkorpal/minimal-gross/blob/main/scripts/CM_Gross/ref/cm_28.txt}{\texttt{cm\_28.txt}}  \\ \hline
  $-96^3$ & $2,3,8,10,12,13,14,15,18\pmod {19}$ & \href{https://github.com/gkorpal/minimal-gross/blob/main/scripts/CM_Gross/ref/cm_19.txt}{\texttt{cm\_19.txt}}  \\ \hline
  $-3\cdot 160^3$ &  $2,5,8,11,14,17,20,23,26\pmod {27}$ & \href{https://github.com/gkorpal/minimal-gross/blob/main/scripts/CM_Gross/ref/cm_27.txt}{\texttt{cm\_27.txt}} \\ \hline
  $-960^3$  & $2, 3, 5, 7, 8, 12, 18, 19, 20, 22, 26, 27, 28, 29, 30,$ $ 32, 33, 34, 37, 39, 42\pmod {43}$ & \href{https://github.com/gkorpal/minimal-gross/blob/main/scripts/CM_Gross/ref/cm_43.txt}{\texttt{cm\_43.txt}}  \\ \hline
  $-5280^3$ & $2, 3, 5, 7, 8, 11, 12, 13, 18, 20, 27, 28, 30, 31, 32,$ $ 34, 38, 41, 42, 43, 44, 45, 46, 48, 50, 51, 52, 53,$ $ 57, 58, 61, 63, 66 \pmod{67}$ & \href{https://github.com/gkorpal/minimal-gross/blob/main/scripts/CM_Gross/ref/cm_67.txt}{\texttt{cm\_67.txt}}  \\ \hline
  $-640320^3$  &  $2, 3, 5, 7, 8, 11, 12, 13, 17, 18, 19, 20, 23, 27, 28, $ $ 29, 30, 31, 32, 37, 42, 44, 45, 48, 50, 52, 59, 63, 66, $ $ 67, 68, 70, 72, 73, 75, 76, 78, 79, 80, 82, 86, 89, 92, $ $ 94, 98, 99, 101, 102, 103, 105, 106, 107, 108, 109, $ $ 110, 112, 114, 116, 117, 120, 122, 123, 124, 125, $ $ 127, 128, 129, 130, 137, 138, 139, 141, 142, 147, $ $ 148, 149, 153, 154, 157, 159, 162\pmod{163}$& \href{https://github.com/gkorpal/minimal-gross/blob/main/scripts/CM_Gross/ref/cm_163.txt}{\texttt{cm\_163.txt}}  \\ \hline
\end{tabular}
\end{table}

%%%%%%%%%%%%%%%%
\bibliographystyle{amsalpha}
\bibliography{refs}
%%%%%%%%%%%%%%%%%%%%%%%%%%%%%%%%%%%%%%%%%%%%%%%%

%%%%%%%%%%%%%%%%%%%%%
\appendix

\section{Proof of \Cref{thm:4types}}\label{appendix}

The authors would like to warmly thank Jonathan Love for providing a much simpler and shorter proof of Theorem \ref{thm:4types}, which we present here.

As in the body of the paper, let  $B_{p}$ be the quaternion algebra over $\mathbb{Q}$ with discriminant $p$, $\mathcal{O}$ be an order in $B_p$ (which will turn out to be maximal), $\mathcal{O}^T$ be the Gross lattice of $\mathcal{O}$, and $\{\beta_1, \beta_2, \beta_3\}\subset \mathcal{O}^T$ be a successive minimal basis for $\mathcal{O}^T$. Furthermore, we write 
\[
D_i := \norm{\beta_i}^2, \qquad \text{and}\qquad t_{ij} := \frac{1}{2}\trd(\beta_i\overline{\beta_j}) = -\frac{1}{2}\trd(\beta_i\beta_j),
\]
where the last equality follows because $\trd(\beta_j) = 0$ for each $j$. As shown in Lemma \ref{lem:Gram_matrix}, we may choose $t_{12}$ and $t_{13}$ to be nonnegative by negating $\beta_2$ or $\beta_3$, as necessary, and in that case we have $0 \leq t_{12}, t_{13} \leq \frac{D_1}{2}$.

\begin{lemma}\label{lem:Omax}
Let $p$ be a prime and $\mathcal{O}$ an order in $B_p$ be such that the lattice $\langle \beta_1,\beta_2\rangle$ spanned by two vectors attaining the first two successive minima $D_1$ and $D_2$ of its Gross lattice $\mathcal{O}^T$ has determinant $4p$. 
    Then the order $\mathcal{O}$ is maximal, and is spanned as a lattice by $\{1, \alpha_1, \alpha_2, \alpha_1\alpha_2\}$, where for $i = 1,2,$ $\alpha_i = \frac{1}{2}(\epsilon_i + \beta_i)$ with $\epsilon_i \in \{0,1\}$ satisfying $\epsilon_i \equiv D_i \pmod{2}$.
\end{lemma}

\begin{proof}
    As argued in \cite[Section 3.3]{Goren-Love}, the elements $\alpha_i$ as constructed here must belong to the order $\mathcal{O}$. As a consequence $\mathcal{O}$ must certainly contain the sublattice $\langle 1, \alpha_1, \alpha_2, \alpha_1\alpha_2\rangle$.
A computation as in \cite[Proposition 3.12]{Goren-Love} shows that the discriminant of the trace matrix of this lattice (by which we mean twice its Gram matrix) is equal to
\[
\frac{(D_1 D_2 - t_{12}^2)^2}{16} = p^2.
\]

This discriminant must be divisible by the discriminant of $\mathcal{O}$, which is itself divisible by the discriminant of a maximal order containing $\mathcal{O}$, which equals $p^2$. Hence all lattices in this chain of inclusions are equal and $\mathcal{O}$ is maximal.
\end{proof}

Now let $B_p^0 = \{x \in B_p : \trd(x) = 0\}$, and given any $v, w \in B_p^0$, define the cross product
\[
v \times w := vw - \frac{1}{2}\trd(vw)
\]
(this is the orthogonal projection of $vw$ onto $B_p^0$). By construction, $v \times w$ has trace $0$ and
norm $\norm{v}^2\norm{w}^2 - \frac{1}{4}\trd(vw)^2$, and it is orthogonal to both $v$ and $w$. Indeed we have
\begin{align*}
\trd(v(v \times w)) & = \trd(v^2w) - \frac{1}{2}\trd(vw)\trd(v)\\
&= -\norm{v}^2\trd(w) - \frac{1}{2}\trd(vw)\trd(v),
\end{align*}
and this vanishes because $\trd(v) = \trd(w) = 0$. A similar computation shows that
$\trd(w(v \times w)) = \trd((v \times w)w) = 0$.

\begin{lemma}\label{lem:3rd_minimum}
Let $p$ be a prime and $\mathcal{O}$ be an order in $B_p$ be such that the lattice $\langle \beta_1,\beta_2\rangle$ spanned by two vectors attaining the first two successive minima $D_1$ and $D_2$ of the Gross lattice $\mathcal{O}^T$ has determinant $4p$. 
    Then the elements $\beta_1, \beta_2,$ and $\beta_3 := \frac{1}{2}(\beta_1 \times \beta_2 + \epsilon_2\beta_1 - \epsilon_1\beta_2)$ are a basis of $\mathcal{O}^T$ attaining its successive minima.
\end{lemma}

\begin{proof}
By Lemma \ref{lem:Omax}, $\mathcal{O}$ is spanned by $\{1, \alpha_1, \alpha_2, \alpha_1\alpha_2\}$ where for $i = 1,2,$ $\alpha_i = \frac{1}{2}(\epsilon_i + \beta_i)$ with $\epsilon_i \in \{0,1\}$ satisfying $\epsilon_i \equiv D_i \pmod{2}$. As a consequence the Gross lattice of $\mathcal{O}$ is spanned by the images of the elements $\alpha_1, \alpha_2,$ and $\alpha_1\alpha_2$ under $x \mapsto 2x - \trd(x)$. The first two of these images are of course $\beta_1$ and $\beta_2$ by construction. 

Now recalling that
\begin{align*}
\trd(\beta_1\beta_2) & = \trd((2\alpha_1 - \trd(\alpha_1))(2\alpha_2 - \trd(\alpha_2))) \\
& = 4\trd(\alpha_1\alpha_2) - 2\trd(\alpha_1)\trd(\alpha_2),
\end{align*}
we have
\begin{align*}
2\alpha_1\alpha_2 - \trd(\alpha_1\alpha_2) & = \frac{1}{2}(\beta_1 + \epsilon_1)(\beta_2 + \epsilon_2) - \trd(\alpha_1\alpha_2)\\
& = \frac{1}{2}(\beta_1\beta_2 + \epsilon_1\beta_2 + \epsilon_2\beta_1 + \epsilon_1\epsilon_2) - \frac{1}{4}\trd(\beta_1\beta_2) - \frac{1}{2}\epsilon_1\epsilon_2\\
& = \frac{1}{2}(\beta_1 \times \beta_2 + \epsilon_1\beta_2 + \epsilon_2\beta_1).
\end{align*}

Now we may subtract $\epsilon_1\beta_2$ from this to obtain $\beta_3$, from which we conclude that $\beta_3 \in \mathcal{O}^T$.

To check that $\beta_3$ attains the third successive minimum of $\mathcal{O}^T$, we must compute the norm
of $\beta_3 + c\beta_1 + d\beta_2$ for $c, d \in \mathbb{Z}$. Writing $x = 2c + \epsilon_2$ and $y = 2d - \epsilon_1$ we have
\begin{align*}
\norm{2(\beta_3 + c\beta_1 + d\beta_2)}^2 & = \norm{\beta_1 \times \beta_2}^2 + \norm{x\beta_1 + y\beta_2}^2 \\
& = 4p + D_1x^2 + 2t_{12}xy + D_2y^2 \\
& = 4p + (D_1 - t_{12})x^2 + t_{12}(x + y)^2 + (D_2 - t_{12})y^2.
\end{align*}

Since $D_1 - t_{12}$, $t_{12},$ and $D_2 - t_{12}$ are all nonnegative, if we can simultaneously minimize
$|x|, |x + y|,$ and $|y|$, then the norm must be minimized.

We consider the following cases:
\begin{itemize}
\item If $\epsilon_1 = \epsilon_2 = 0$ (so $x, y \in 2\mathbb{Z}$), $|x|, |x + y|, |y|$ are simultaneously minimized when
$x = y = 0$.

\item If $\epsilon_1 = 0$ and $\epsilon_2 = 1$ (so $y \in 2\mathbb{Z}$ but $x \in 2\mathbb{Z} + 1$) we will necessarily have $|x|, |x + y| \geq 1$, so the minimum is attained for $x = 1$ and $y = 0$.

\item If $\epsilon_1 = 1$ and $\epsilon_2 = 0$, we have $|y|, |x + y| \geq 1$; so $y = -1$ and $x = 0$ attains the minimum. 

\item If $\epsilon_1 = \epsilon_2 = 1$, then $|x|, |y| \geq 1$ and $|x + y| \geq 0$; the minimum is attained with $x = 1$ and $y = -1$.
\end{itemize}

Therefore we see that in all cases the minimum is attained for $c = d = 0$, so $\beta_3$ attains the third successive
minimum of $\mathcal{O}^T$. 
    
\end{proof}

We are now in a position to give a proof of Theorem \ref{thm:4types}:

\begin{proof}[Proof of \Cref{thm:4types}]
With the notation established so far in this appendix, we have:
\begin{gather*}
\trd(\beta_1\bar{\beta}_3) = -\frac{1}{2}\trd(\epsilon_2\beta_1^2 - \epsilon_1\beta_1\beta_2)
= \epsilon_2D_1 - \epsilon_1 t_{12}, \\
\trd(\beta_2 \bar{\beta}_3)
= -\frac{1}{2}\trd(\epsilon_2 \beta_2 \beta_1 - \epsilon_1 \beta_2^2)
= \epsilon_2 t_{12} - \epsilon_1 D_2,\\
\norm{\beta_3}^2
= p + \frac{1}{4}(D_1 \epsilon_2^2 - 2 t_{12} \epsilon_1 \epsilon_2 + D_2 \epsilon_1^2).
\end{gather*}

Recalling that $\epsilon_i \in \{0,1\}$ so $\epsilon_i^2 = \epsilon_i$, a Gram matrix for $\mathcal{O}^T$ is thus
\begin{equation*}
\begin{pmatrix}
D_1 & t_{12} & \frac{1}{2}(\epsilon_2 D_1 - \epsilon_1 t_{12}) \\
t_{12} & D_2 & \frac{1}{2}(\epsilon_2 t_{12} - \epsilon_1 D_2) \\
\frac{1}{2}(\epsilon_2 D_1 - \epsilon_1 t_{12}) & \frac{1}{2}(\epsilon_2 t_{12} - \epsilon_1 D_2) & p + \frac{1}{4}(D_1 \epsilon_2 - 2 t_{12} \epsilon_1 \epsilon_2 + D_2 \epsilon_1)
\end{pmatrix}.
\end{equation*}

Now we have by hypothesis that $t_{12} = \frac{1}{2}\trd(\beta_1 \bar{\beta}_2) = -\frac{1}{2}\trd(\beta_1 \beta_2)$ satisfies $D_1 D_2 - t_{12}^2 = 4p$ and hence $t_{12} \equiv D_1 D_2 \pmod{2}$; recall that we also have that $0 \leq t_{12} \leq \frac{D_1}{2}$, $\epsilon_i \equiv D_i \pmod{2}$, and $D_i \equiv 0, 3 \pmod{4}$. This naturally gives rise to the following four possibilities:

\begin{itemize}

\item \textbf{Type 1:} $D_1 \equiv D_2 \equiv 0 \pmod 4$. Then $\epsilon_1 = \epsilon_2 = 0$, and the normalized Gram matrix for $\mathcal{O}^T$ is
\[
\begin{pmatrix}
D_1 & t_{12} & 0 \\
t_{12} & D_2 & 0 \\
0 & 0 & p
\end{pmatrix}.
\]

\item \textbf{Type 2:} $D_1 \equiv 0 \pmod{4}$ and $D_2 \equiv 3 \pmod 4$. Then $\epsilon_1 = 0$ and $\epsilon_2 = 1$, and the normalized Gram matrix for $\mathcal{O}^T$ is
\[
\begin{pmatrix}
D_1 & t_{12} & \tfrac{D_1}{2} \\
t_{12} & D_2 & \tfrac{t_{12}}{2} \\
\tfrac{D_1}{2} & \tfrac{t_{12}}{2} & p + \tfrac{D_1}{4}
\end{pmatrix}.
\]

\item \textbf{Type 3:} $D_1 \equiv 3$ and $D_2 \equiv 0 \bmod 4$. Then $\epsilon_1 = 1$ and $\epsilon_2 = 0$, and a Gram matrix for $\mathcal{O}^T$ is
\[
\begin{pmatrix}
D_1 & t_{12} & -\tfrac{t_{12}}{2} \\
t_{12} & D_2 & -\tfrac{D_2}{2} \\
-\tfrac{t_{12}}{2} & -\tfrac{D_2}{2} & p + \tfrac{D_2}{4}
\end{pmatrix}.
\]

To obtain the normalized Gram matrix, we can replace $\beta_3$ with $-\beta_3$ to get a Gram matrix where the top-right entry is nonnegative:
\[
\begin{pmatrix}
D_1 & t_{12} & \tfrac{t_{12}}{2} \\
t_{12} & D_2 & \tfrac{D_2}{2} \\
\tfrac{t_{12}}{2} & \tfrac{D_2}{2} & p + \tfrac{D_2}{4}
\end{pmatrix}.
\]

\item \textbf{Type 4:} $D_1 \equiv D_2 \equiv 3 \bmod 4$. Then $\epsilon_1 = \epsilon_2 = 1$, and the normalized Gram matrix for $\mathcal{O}^T$ is
\[
\begin{pmatrix}
D_1 & t_{12} & \tfrac{D_1 - t_{12}}{2} \\
t_{12} & D_2 & \tfrac{t_{12} - D_2}{2} \\
\tfrac{D_1 - t_{12}}{2} & \tfrac{t_{12} - D_2}{2} & p + \tfrac{D_1 - 2t_{12} + D_2}{4}
\end{pmatrix}.
\]

Note that $t_{12} \equiv D_1 D_2 \equiv 1 \pmod 2$, so $D_1 - t_{12}$ and $t_{12} - D_2$ are both even, and $D_1 + D_2 - 2t_{12} \equiv 0 \bmod 4$, and therefore again all entries of this matrix are integers.
\end{itemize}
   
\end{proof}
\end{document}